\documentstyle[12pt,amstex,epsf]{article}

\catcode`\@=11
\def\numberbysection{\@addtoreset{equation}{section}
        \def\theequation{\thesection.\arabic{equation}}}
\numberbysection
\newtheorem{theo}{Theorem}
\newtheorem{lemma}{Lemma}
\newtheorem{prop}{Proposition}

\def\bea{\begin{eqnarray}}
\def\eea{\end{eqnarray}}
\def\choice{{\vec{\varepsilon}}}
\def\orga{{{\mathrm or}\,}_\Gamma}
\def\hz{\hat{Z}}
\def\R{{\Bbb{R}}}
\def\shuf{\times_\sigma}
\begin{document}
\begin{titlepage}
\begin{center}
\vskip 70pt
{\large \bf Vassiliev knot invariants
and Chern-Simons perturbation theory to all orders}
\vskip 0.5in
Daniel Altschuler
\footnote{Supported by Fonds national suisse de la
recherche scientifique.}\\
{\em Institut f\"ur Theoretische Physik\\
ETH-H\"onggerberg\\
CH-8093 Z\"urich, Switzerland}
\vskip 0.2in
and 
\vskip 0.2in
Laurent Freidel\\
{\em Laboratoire de Physique Th\'eorique ENSLAPP
\footnote{URA 1436 du CNRS, associ\'ee \`a l'Ecole Normale
Sup\'erieure de Lyon et \`a l'Universit\'e de Savoie.}\\
Ecole Normale Sup\'erieure de Lyon\\
46, all\'ee d'Italie,
69364 Lyon Cedex 07, France}
\end{center}
\vskip .5in
\begin{abstract}
At any order, the perturbative expansion of the 
expectation values of Wilson lines in Chern-Simons theory 
gives certain integral expressions.
We show that they all lead 
to knot invariants. Moreover these
are finite type invariants whose order
coincides with the order in the perturbative expansion.
Together they combine to give a universal Vassiliev invariant.
\end{abstract}
\vskip .6in
\begin{center}
Revised version, October 1996, to appear in Commun. Math. Phys.
\end{center}
\end{titlepage}
\pagenumbering{arabic}
\section{Introduction}
Chern-Simons theory is the most popular example of topological
field theory in 3 dimensions. Given 
a compact Lie group $G$, a compact, oriented 3-manifold
$M$, a link $L\subset M$, and for each component of $L$ a 
representation of $G$, this theory associates
topological invariants to these data. There are several ways
to define the invariants, which are all closely related.
First of all there are the non-perturbative definitions:
Witten \cite{witten} used fundamental properties of quantum field theory,
in particular the path integral formulation, and Reshetikhin and Turaev
\cite{rt} used quantum groups. These two definitions are equivalent.

Then there are the perturbative definitions,
the first of which were given by Guadagnini et al. \cite{guad}
in the case $M=S^3$, $L\neq\emptyset$, using propagators and
Feynman diagrams. This approach was then elaborated by
Bar-Natan \cite{bn1,bn11}. The case $M\neq S^3$, $L=\emptyset$ was treated
by Axelrod and Singer \cite{as}. A common feature of all these
works is the Feynman diagram expansion familiar in
perturbative quantum field theory. Invariants are defined
at every order in the expansion, each is a sum of
several terms corresponding to the diagrams of the given order.
The contribution of any
diagram is the product of two factors, the first 
depends only on the group $G$ and the representations associated
to the components of $L$, and the second is independent of
$G$ and its representations, it is an integral over the configuration
space of the vertices of the diagram, some of which are constrained to
lie on $L$, while the others can lie anywhere in the complement of $L$.
When $L$ is a knot in $S^3$, several properties of
the invariant arising from the contributions
of order two were already discussed in \cite{guad}, although
the invariance itself was shown in \cite{bn11}.
Bar-Natan also studied the properties
of the group-dependent contributions, and among them he
found relations between the contributions of different diagrams
which are the same for all groups $G$. This led him \cite{bn2} to define
abstract objects, which we call BN diagrams, by
these relations, and abstract invariants
which take their values in the space of BN diagrams.
To every choice of group $G$ and representations corresponds
a linear functional on the space of BN diagrams. Applying
this functional to the abstract invariants gives back
the ordinary group-dependent invariants.

In order to show that the contributions of a given order 
sum up to an invariant, one must compute the variation of these
integrals under a small change of the embedding of $L$, and this 
proved to be quite difficult and lengthy. However,
Bott and Taubes \cite{bt} greatly improved this situation.
They showed that the variation can be split in two terms,
the ``diagrammatic'' and the ``anomalous'' variations.
As its name indicates, the diagrammatic variation can be read
at once from the Feynman diagram. It corresponds to the
differential of Kontsevich's graph complex, obtained by
collapsing the edges. The anomalous variation is more
difficult to compute, but it is proportional
to the variation of the first order contribution, the
``self-linking number''. The constant of proportionality,
is still unknown in general, but independent
of the embedding. These results of Bott and Taubes are powerful
enough, as we will show, to prove invariance at all orders.

During the same period, the subject of Vassiliev knot invariants,
also known as finite type invariants, was developing rapidly. 
The starting point of Vassiliev \cite{vassil} was 
the space of all immersions of $S^1$ in $S^3$. In this space,
a knot type is a cell whose faces are singular knots with
a finite number of transverse double points. Any knot
invariant can be extended to such singular knots. It is said to
be a finite type invariant of order $\leq N$, if it vanishes on
all singular knots with more than $N$ double points.
Let $V^N$ be the space of invariants of order $\leq N$.
Unexpectedly at first, Bar-Natan found that $V^N/V^{N-1}$ 
embeds in the dual of the space of BN diagrams of degree $N$. 
Kontsevich \cite{konts} showed that the two spaces are in fact isomorphic.
His proof \cite{bn2} involved the construction of a universal Vassiliev
invariant, a formal power series in the space of BN diagrams
whose coefficients are finite type invariants, based on the
Knizhnik-Zamolodchikov equations of the WZW model of
conformal field theory. (It was Witten \cite{witten} who discovered the
relation between conformal field theory and topological field
theory. Fr\"ohlich and King \cite{fk} were the first to construct
link invariants from the KZ equations.)

In this paper we start from the results of Bott and Taubes \cite{bt} to 
construct a universal Vassiliev knot invariant, given by the
perturbative expansion of the expectation value of a Wilson loop in 
Chern-Simons theory on $\R^3$. The basic ingredient in the integrals
obtained from the Feynman rules is the
propagator of the gauge field, which is given in the Lorentz gauge
by the Gauss two-form, the pullback of the volume form on $S^2$.
During the Aarhus conference on geometry and physics in the summer
of 1995, we learned from Dylan Thurston that he independently 
obtained similar results to ours \cite{thurston}. Recently,
perturbative Chern-Simons invariants have also been investigated
in \cite{alvlab,hisa}.

In more details, the contents of the paper are as follows:
in section 2, we define the graphs appearing in the perturbative expansion,
which are equipped with an additional structure called vertex
orientation, and state some simple combinatorial lemmas. In section
3 we give the Feynman rules, in which the vertex orientation plays
an important role. They allow us to define unambiguously the signs
of the contributions of graphs appearing in the perturbative expansion
(see propositions \ref{promor} and \ref{proepsi}). In section 4 we
define the expectation value of a Wilson loop $Z$, which is a sum over
trivalent graphs, and prove that it is invariant under the changes
of embedding corresponding to the collapse of a single edge
(theorem \ref{delta}). In section 5 we consider the other variations
of the embedding, called ``anomalous''. We improve some results
of \cite{bt} in proposition \ref{proepsi} and theorem \ref{annul},
which allow us to conclude in theorem \ref{invariant}
that a suitably corrected version of
$Z$ becomes a framed knot invariant $\hz$. In the last section
we prove that $\hz$ is a universal invariant. In particular, the
$N$\/-th order contribution to $\hz$ is a 
finite type invariant of order $\leq N$. Although it is stated
explicitly in the literature, we have never seen a proof of this
going beyond the second order. 
The question whether the KZ and the Chern-Simons 
universal invariants are equal is still open. The answer would be
positive if one could show that the Chern-Simons invariant 
extends functorially to the category of tangles, as in the case
of the KZ invariant, but at least to us it is not obvious that it
has this extension property. In an appendix, we recall the definition
of the pushforward, or integration along the fiber, which enters
the formulation of the Feynman rules.

\section{Graphs}
A Wilson graph $\Gamma$ is a one-dimensional, connected, simplicial
complex equipped with some additionals structures, which we now
describe.  We assume that all vertices have valence $\geq 3$ and that
each graph has a distinguished oriented cycle $W_{\Gamma}$, called the
Wilson line.  We define:
\begin{equation}
\begin{array}{ll}
V_{\Gamma} = \left\{ \mbox{vertices of ${\Gamma}$} \right\}, &
E_{\Gamma} = \left\{\mbox{edges of ${\Gamma}$ }\right\}, \\
V_{\Gamma}^{o} = V_{\Gamma} \bigcap W_{\Gamma}, &
E_{\Gamma}^{o} = \left\{ e \in E_{\Gamma} | e \not \subset  W_{\Gamma},
\partial e = e \bigcap W_{\Gamma}\right\},  \\
V_{\Gamma}^{i} = V_{\Gamma} - V_{\Gamma}^{o}
, & 
E_{\Gamma}^{a} = E_{\Gamma} - E_{\Gamma}^{o}, \\
E_{\Gamma}^{W} = E_{\Gamma} \bigcap W_{\Gamma}, &
E_{\Gamma}^{nW} = E_{\Gamma} - E_{\Gamma}^{W} .\\
\end{array}
\end{equation}
$E_{\Gamma}^{a}$ is called the set of admissible edges of $\Gamma$,
$V_{\Gamma}^{o}$ the set of external vertices, 
$V_{\Gamma}^{i}$ the set of internal vertices.
An example of Wilson graph is given in fig. 1. 
The edges of $E_{\Gamma}^{W}$
are solid lines, those of $E_{\Gamma}^{nW}$ are dashed lines.
\begin{figure}[htbp]
$$\epsf{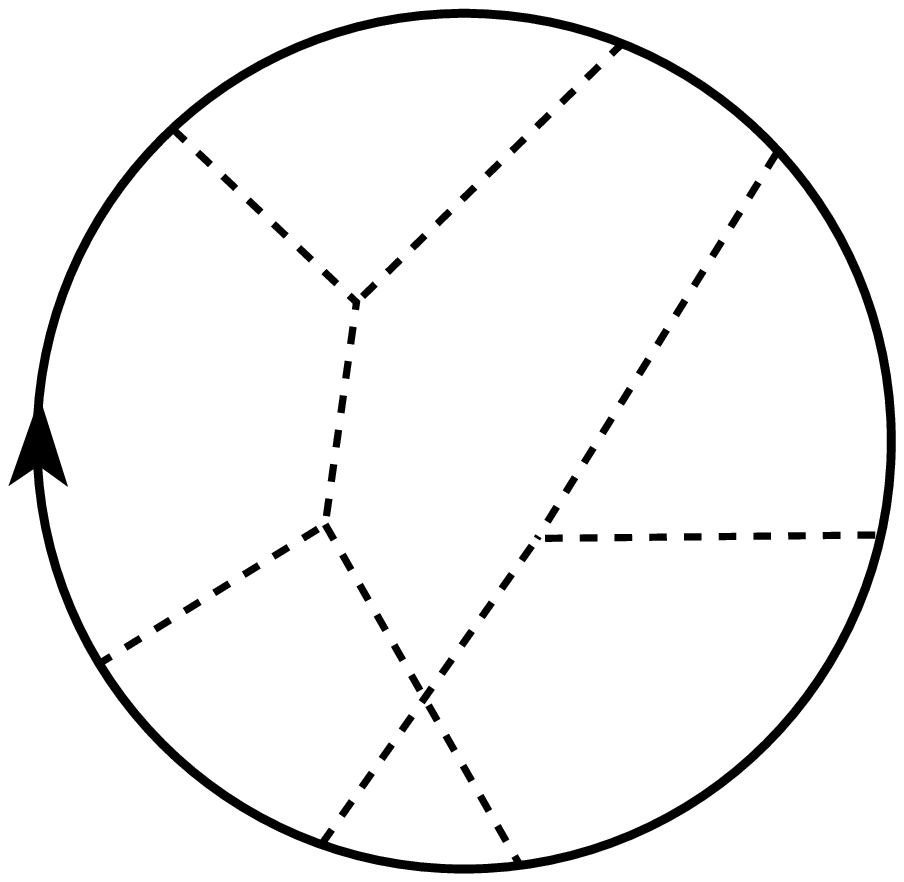}[xscale=1/2,yscale=1/2] 
$$ \caption{\label{wilson}}
\end{figure}
We say that two Wilson graphs are equivalent, and write 
$\Gamma \sim \Gamma^\prime$ , if there are two bijections 
\begin{equation}
\rho : V_{\Gamma} \longrightarrow V_{\Gamma^\prime},\, \sigma : E_{\Gamma} 
\longrightarrow E_{\Gamma^\prime},
\end{equation}
such that  $\sigma|_{E_{\Gamma}^{W}} $ is a map from $E_{\Gamma}^{W}$ to 
$E_{\Gamma^\prime}^{W}$ preserving the orientation of the Wilson line,
and 
\begin{equation}
\rho \,\partial_{\Gamma} =\partial_{\Gamma^\prime}\, \sigma,
\end{equation}
where $\partial_{\Gamma}$ is the usual boundary map acting on 
$E_{\Gamma}$. The pair $(\rho,\sigma)$ is an isomorphism of graphs.
It is an automorphism if $\Gamma=\Gamma'$. Let ${\mathrm Aut}(\Gamma)$
be the finite group of automorphisms of $\Gamma$.

We define an additional structure on graphs which we call vertex
orientation, or V-orientation, which should not be confused with the
usual definition of graph orientation.  We orient $\Gamma$ by 
giving for each vertex an order
among the edges arriving at this vertex. More precisely for each vertex $v
\in V_\Gamma$ we choose a  bijection :
\begin{equation}
o_v : E_v \longrightarrow \{1,2\cdots , |E_v|\},
\label{ov}
\end{equation}
where $E_v=\{ e\in E_\Gamma^{nW}\,|\, v\in \partial e \}$.
We say that $O_\Gamma=\{(o_v)_{v\in V_\Gamma }\}$ and
$O^{\prime}_{\Gamma}=\{(o^{\prime}_v)_{v\in V_\Gamma }\}$ 
define the same 
(resp. opposite) V-orientation if 
\begin{equation}
\prod_{v\in V_\Gamma }\mbox{sign}((o_v)^{-1}o^\prime_v)=+1 \,(\mbox{resp.} -1),
\end{equation}
where sign is the signature of a permutation.
For short, in the sequel $\Gamma$ will denote the graph 
$\Gamma$ with a choice of V-orientation $O_\Gamma$, and $-\Gamma$
the graph equipped with the opposite V-orientation.

Using our definition of isomorphism of graphs, we can define the notion of 
induced orientation:
$(\rho,\sigma) O_\Gamma =
\{(o_v \circ \sigma)_{v\in V_\Gamma }\}$.
Since $\Gamma$ is determined by the couple 
$(\partial_\Gamma, O_\Gamma)$ we put
\begin{equation}
(\rho,\sigma ) \cdot \Gamma = (\rho \partial_\Gamma \sigma^{-1},(\rho,\sigma) 
O_\Gamma ).
\end{equation}
The equivalence relation naturally extends to these couples.

\begin{figure}[htbp]
$$\epsf{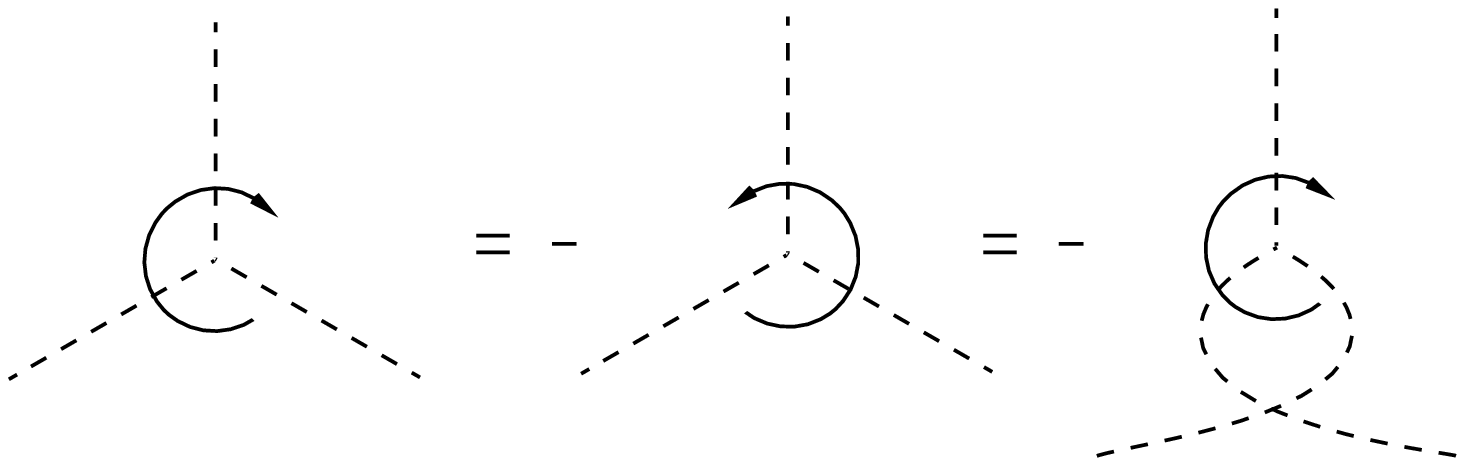}[xscale=2/3,yscale=2/3] 
$$ \caption{\label{orient3}}
\end{figure}

\begin{figure}[htbp]
$$\epsf{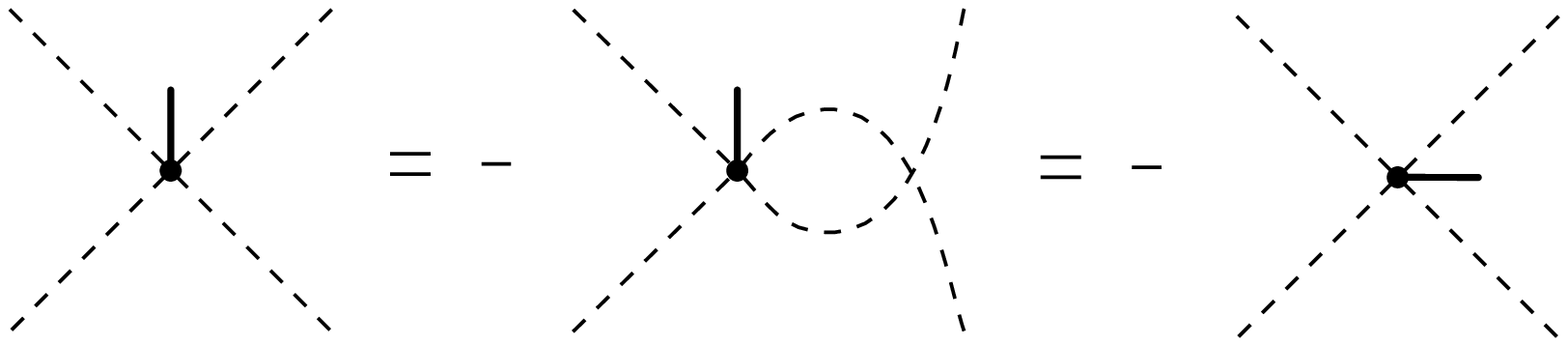}[xscale=2/3,yscale=2/3] 
$$ \caption{\label{orient4}}
\end{figure}

\begin{figure}[htbp]
$$\epsf{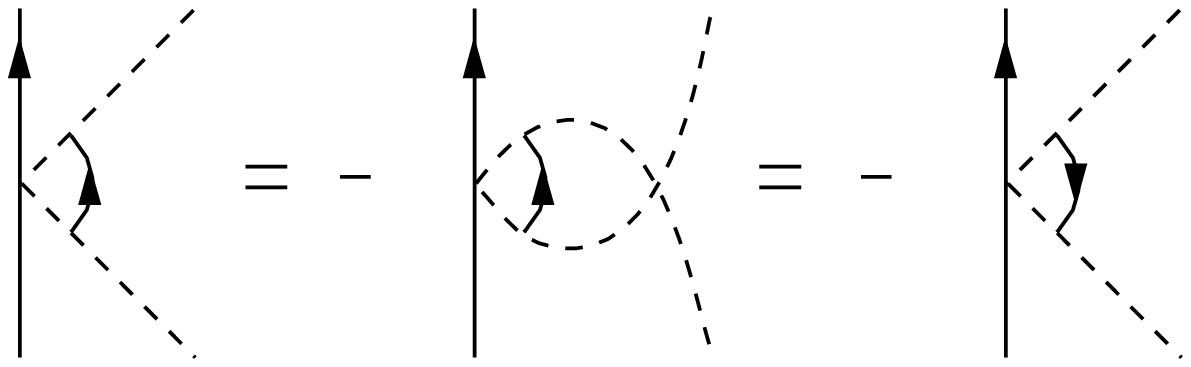}[xscale=2/3,yscale=2/3] 
$$ \caption{\label{orient42}}
\end{figure}

Thus we can say that an automorphism preserves the orientation if the 
induced orientation coincides with the original one, and denote by 
$\mbox{Aut}_+ (\Gamma) $ the normal subgroup of orientation preserving
automorphisms.
In the sequel we will deal only with trivalent and tetravalent vertices. 
In this case it is possible to give a graphical representation of 
V-orientation.
For trivalent vertices the orientation is given by a cyclic order, 
see fig. \ref{orient3}.
For internal tetravalent vertices, 
because the signature of a cyclic permutation 
of 4 objects is $-1$, the orientation 
is given by a separation between the first and the last 
edges, see fig. \ref{orient4}.
For Wilson tetravalent vertices,  
the orientation is given by an order between the internal edges,
see fig. \ref{orient42}.
We will not draw the cyclic orientation of trivalent vertices if they all 
agree with the orientation of the Wilson loop.

\begin{figure}[htbp]
$$\epsf{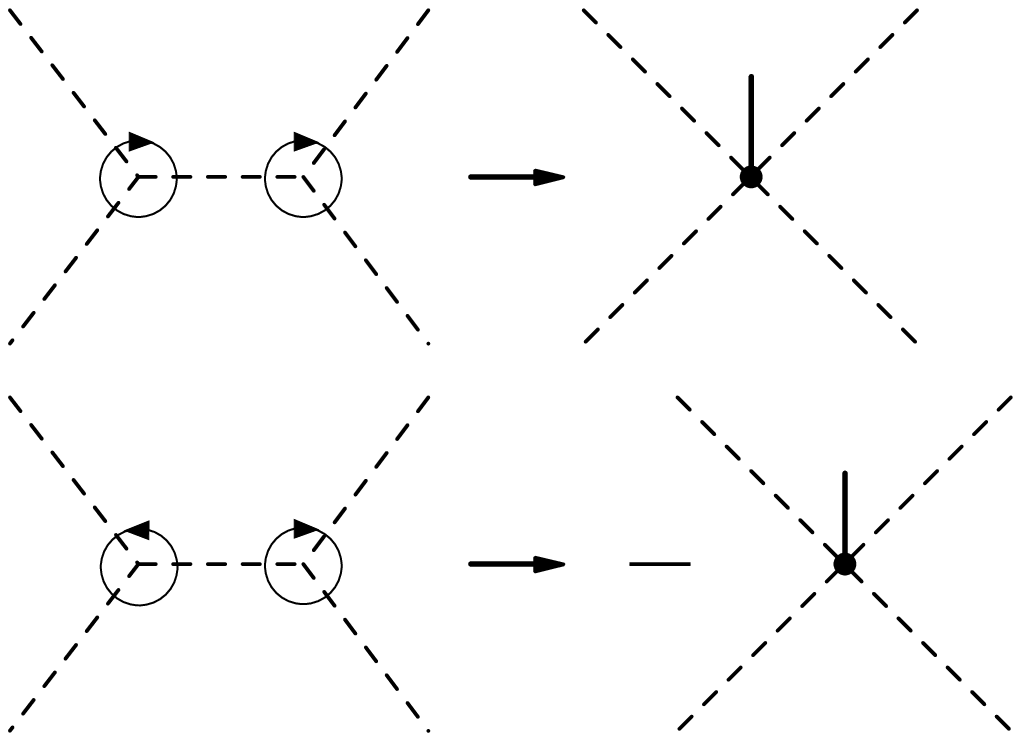}[xscale=2/3,yscale=2/3] 
$$ \caption{\label{contrac1}}
\end{figure}

\begin{figure}[htbp]
$$\epsf{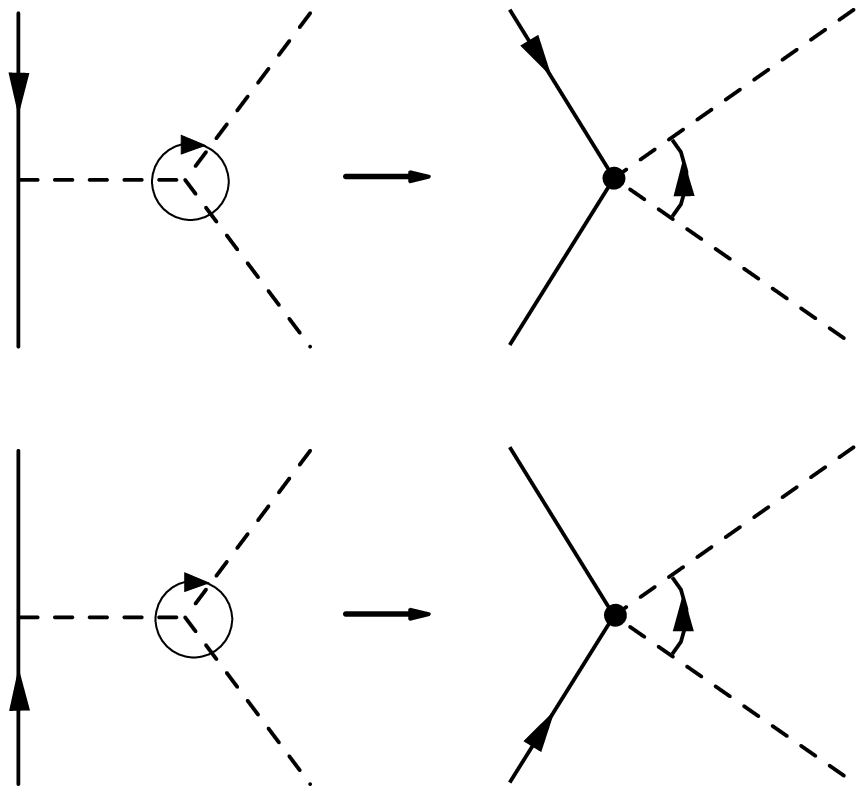}[xscale=2/3,yscale=2/3] 
$$ \caption{\label{contrac2}}
\end{figure}

Let $\Gamma$ be a trivalent graph, and $e \in E_{\Gamma}^{a}$
an admissible edge. We define $\delta_e \Gamma$ to be the graph
obtained from $\Gamma$ by collapsing $e$ to a 4-valent vertex $x$, and 
the vertex 
orientation of $\delta_e \Gamma$ is defined by the following rule:
 
If $e \bigcap W_\Gamma = \emptyset$ we define it by fig. \ref{contrac1}.
If  $e \bigcap W_\Gamma \not= \emptyset$ and 
$e \not\subset W_\Gamma $ we define it by fig. \ref{contrac2}.

\begin{figure}[htbp]
$$\epsf{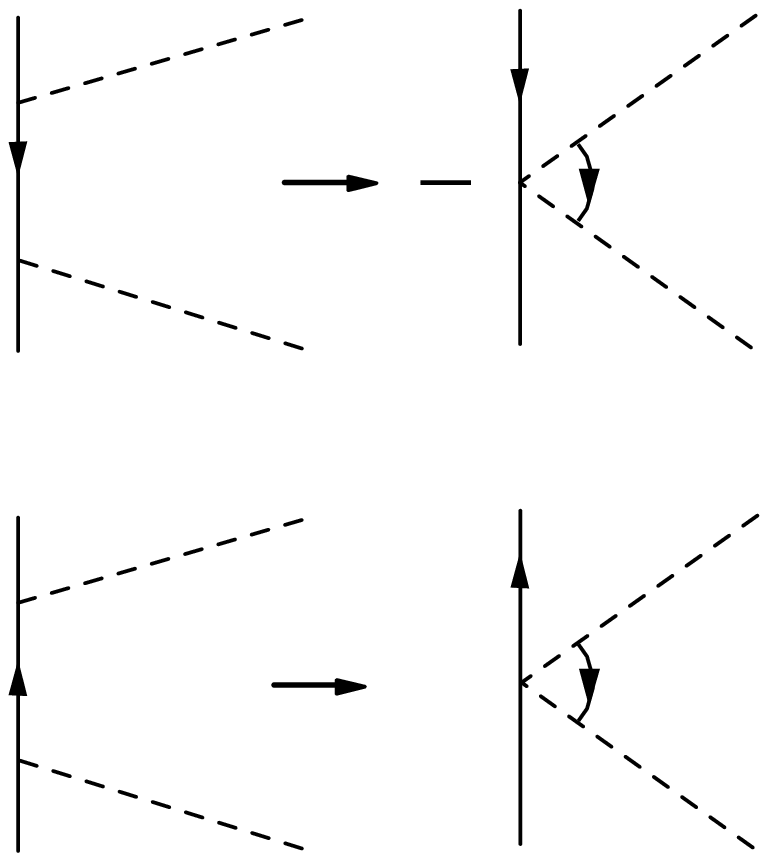}[xscale=2/3,yscale=2/3] 
$$ \caption{\label{contrac3}}
\end{figure}

\begin{figure}[htbp]
$$\epsf{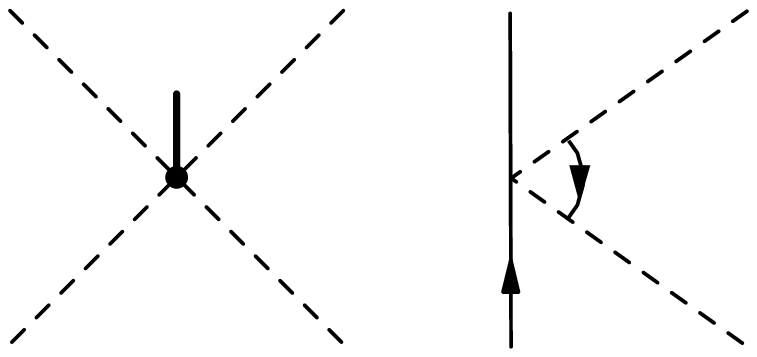}[xscale=2/3,yscale=2/3] 
$$ \caption{\label{quadri}}
\end{figure}

\begin{figure}[htbp]
$$\epsf{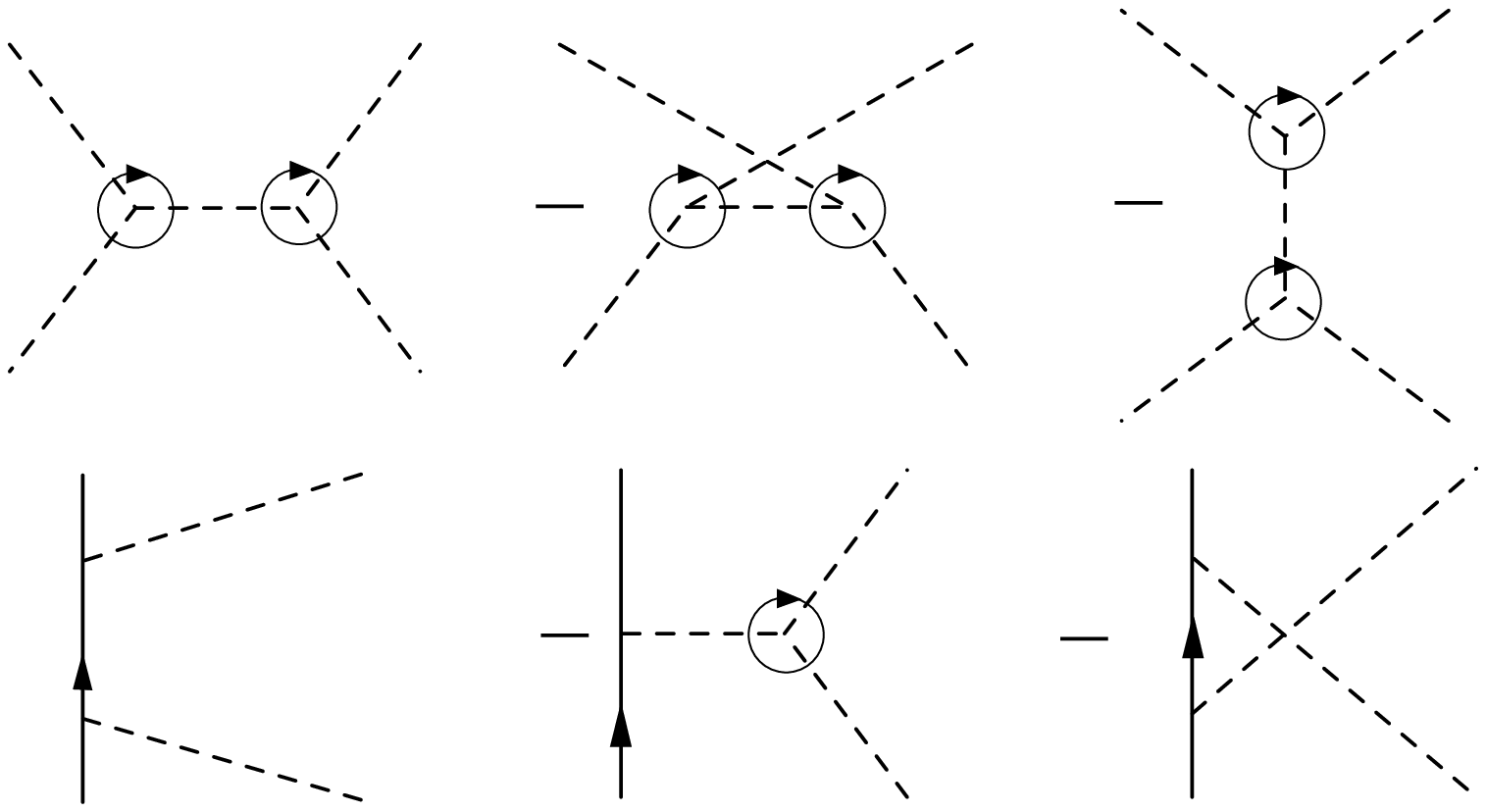}[xscale=2/3,yscale=2/3] 
$$ \caption{\label{relevement}}
\end{figure}

If $e \subset  W_\Gamma $ we define it by fig. \ref{contrac3},
where the oriented edges are those of $W_\Gamma$.
\begin{lemma}
Let  $ \Gamma_x $ be a graph with one 4-valent vertex x, and all other 
vertices trivalent. Then there are at most three 3-valent graphs 
$\Gamma_{x}^i,\, i=1,2,3$, and edges $e_i \in E_{\Gamma^i_x}$ such that 
$\delta_{e_i} \Gamma_{x}^i =\Gamma_{x}$.
If the vertex x of $\Gamma_{x}$ looks like fig. \ref{quadri},
then the vicinity of the edges $e_i$ in $\Gamma_{x}^i$ is given by
fig. \ref{relevement}.
\label{croix}
\end{lemma}
However, it may happen that $\Gamma_{x}^i\sim \Gamma_{x}^j$ for $i\neq j$.
\begin{lemma} 
Let $\Gamma$ a trivalent graph  $e,f \in E_{\Gamma}^a$. 
Then $\delta_{e} \Gamma \sim \delta_{f} \Gamma$ iff there
exists $(\rho,\sigma)\in {\mathrm Aut}(V_\Gamma) \times 
{\mathrm Aut}(E_\Gamma)$ such that
\bea
\delta_e ((\rho,\sigma)\cdot\Gamma) & = & \delta_e \Gamma, \\
\sigma (f) & = & e. \nonumber
\eea
\label{lemef}
\end{lemma}
We omit the proofs of these lemmas, which are essentially obvious.
Later we shall use the following combination of the two lemmas.
\begin{prop}
Let $\Gamma$ be a trivalent graph, $e,f \in E_{\Gamma}^a$.
Then $\delta_{e} \Gamma \sim \pm \delta_{f} \Gamma$ iff there exists 
$i \in \{1,2,3\}$ and 
$(\rho,\sigma)\in {\mathrm Aut}(V_\Gamma) \times 
{\mathrm Aut}(E_\Gamma)$ such that
\begin{equation}
(\rho,\sigma ) \cdot \Gamma = \pm \Gamma_{x}^i, \quad
{\mathrm and}\ \sigma(f) = e, \ {\mathrm where}\ 
\Gamma_{x}=\delta_e \Gamma .
\end{equation}
\end{prop}

It is easy to see that if $\Gamma$ is a trivalent graph with $n$
external vertices and $t$ internal vertices, $n+t$ is even.
The degree (order) of $\Gamma$ is defined to be:
\begin{equation}
{\mathrm deg}\Gamma = (n+t)/2.
\end{equation}
Let $G^3$ be the set of equivalence classes of vertex-oriented
trivalent graphs.
Consider the vector space $\tilde{\cal{A}}$ over $\R$ 
with basis elements $ \Gamma \in G^3$, and $\cal{A}$ the
quotient of $\tilde{\cal{A}}$ by the subspace spanned by all the vectors 
\begin{equation}
\begin{array}{cl}
(*)  &  \sum_{i=1}^3  \Gamma_{x}^i \\
(**)  &  (\Gamma) + (-\Gamma),
\end{array}
\end{equation}
where $ \Gamma_{x}$ is any graph with a single tetravalent vertex.
We call ${\cal A}$ the space of BN diagrams.
Denote by $D$ the canonical projection from $\tilde{\cal{A}}$ to  
$\cal{A}$.
Note that $(*)$  coincides with the IHX and STU relations of 
Bar-Natan \cite{bn1,bn2}.

\begin{figure}[htbp]
$$\epsf{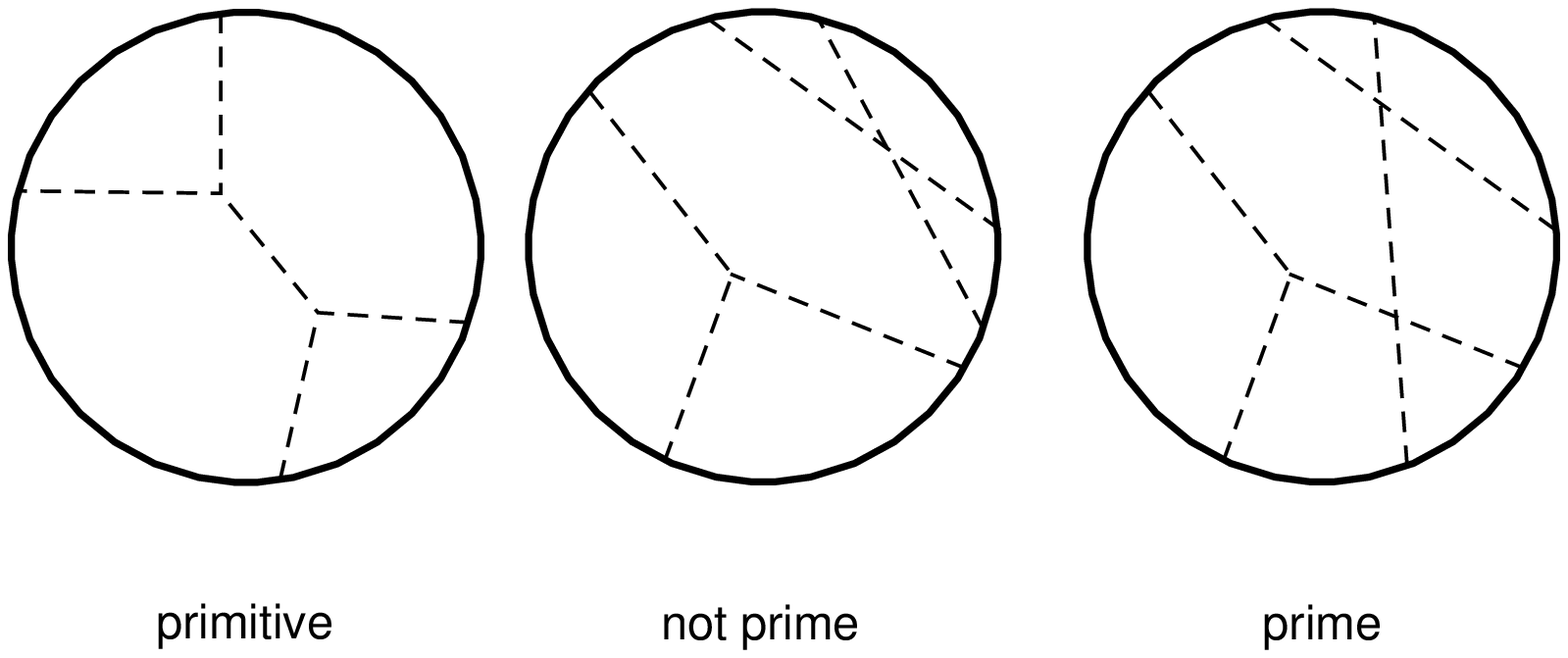}[xscale=2/3,yscale=2/3] 
$$ \caption{\label{prime}}
\end{figure}

It was shown in \cite{bn2} that $\cal{A}$ is a commutative
associative algebra, with the product obtained by taking the
connected sum along the Wilson lines of two graphs. A $prime$
Wilson graph is a graph which cannot be expressed as a product of two 
non-trivial graphs. 
Equivalently, a prime graph is such that if we cut it along any two edges of 
the Wilson line, the resulting graph is connected.
A graph $\Gamma$ is $primitive$ if $\Gamma - W_{\Gamma} $ is connected.
Of course, a primitive graph is prime (see examples in fig. \ref{prime}).
Later we will need a lemma due to Bar-Natan \cite{bn1}:
\begin{lemma}
There exists a unique linear endomorphism $C$
on the space $\cal{A}$ of Bar-Natan diagrams such that:
\begin{eqnarray}
C(D(\Gamma)) & = & D(\Gamma) 
\mbox{  if  $\Gamma$  is a primitive graph, } \\
C(D(\Gamma)) & = & 0 \mbox{ if $\Gamma$ is not prime.}
\label{C}
\end{eqnarray}
\label{lemC}
\end{lemma}
{\em Proof.}
This follows from the fact that $\cal A$ can also be equipped with 
a cocommutative coproduct, so that it becomes a Hopf algebra, and
the primitive graphs correspond to the primitive elements of $\cal A$. 
This implies that $\cal A$ is isomorphic to the enveloping algebra
of the Lie algebra of primitive elements.
Thus the conditions of the lemma determine the value of
$C$ on the elements $\prod_i D(\Gamma_i)$, $\Gamma_i$ primitive,
spanning $\cal A$. $\Box$

\section{Feynman rules}
Let $\cal{K}$ be the space of
all smooth embeddings of $S^1 $ into $\R^3$. 
To any vertex-oriented Wilson graph $\Gamma$  we can associate a
differential form on the space $\cal{K}$, denoted by $I(\Gamma)$. The degree
of this form is equal to:
\begin{equation}
\mbox{degree } I(\Gamma) = \sum_{v \in V_\Gamma} (k_v-3),
\end{equation}
where $k_v$ is the valence of the vertex $v$. Thus for 
trivalent graphs, $I(\Gamma)$ is just a function on the space
of knots $\cal{K}$.
An important property of $I(\Gamma)$ is its behaviour under a change of
vertex orientation:
\begin{equation}
I(-\Gamma)=-I(\Gamma).
\label{reverse}
\end{equation}

\begin{figure}[htbp]
$$\epsf{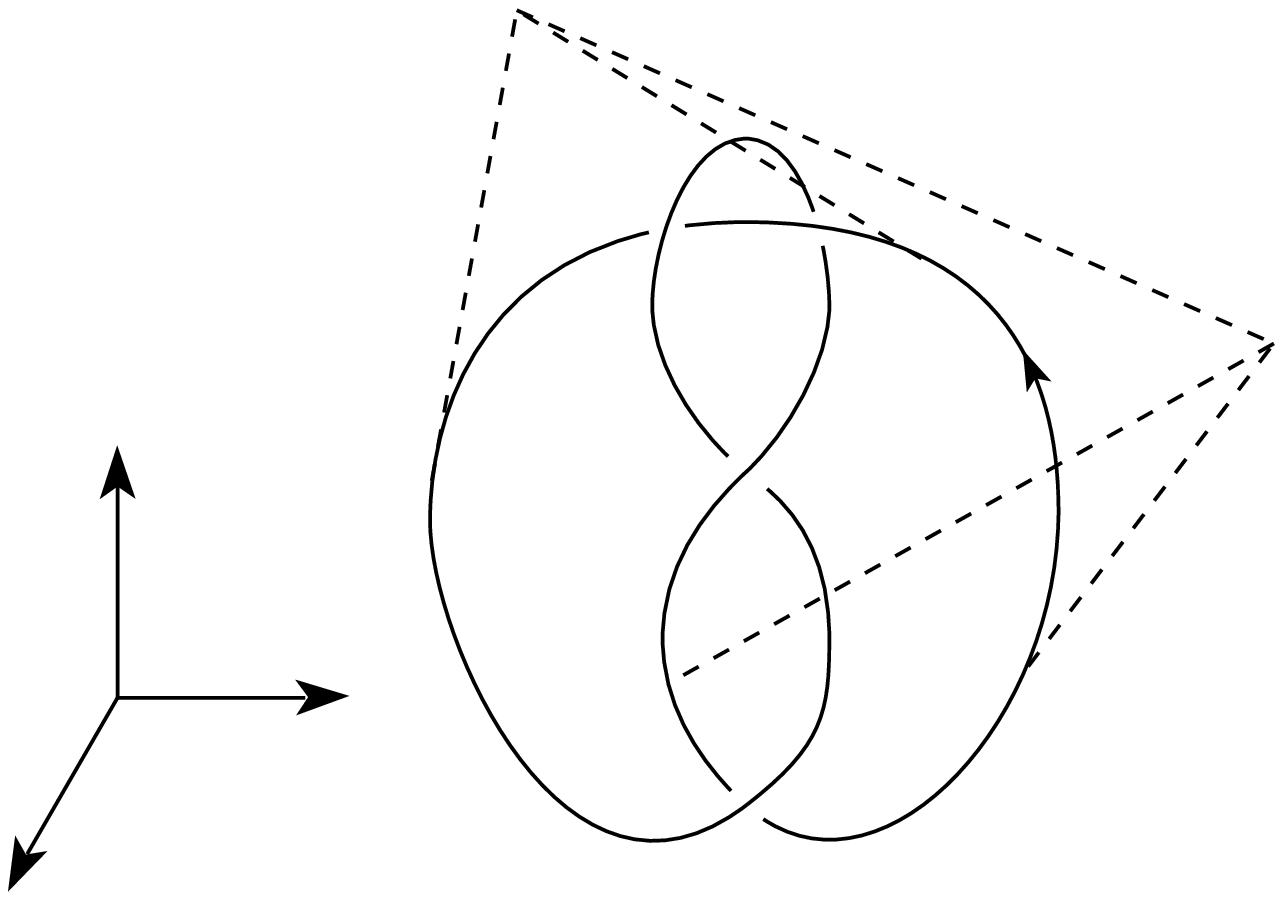}[xscale=2/3,yscale=2/3] 
$$ \caption{\label{huit}}
\end{figure}

The main point of the construction of $I(\Gamma)$ is due to Bott and Taubes.
Let us briefly recall their construction.
First, to any Wilson graph $\Gamma$ we can associate a space $C_\Gamma$,
which is a fiber bundle over the space $\cal{K}$, whose fiber is a compact 
manifold with corners of dimension $n+3t$, where $n$ is the number of Wilson 
vertices of $\Gamma$ and $t$ is the number of internal vertices:
\begin{equation}
p_\Gamma : C_\Gamma \longrightarrow {\cal K}.
\end{equation}
The fiber $p_\Gamma^{-1} (\phi) =\bar{C}_{n,t}(\phi)$, 
where 
$\phi: S^{1} \rightarrow  \R^{3} \, \in \cal{K}$
is a knot,
is a compactification 
of the configuration space $C_{n,t}(\phi)$ of $n$ points on the knot and $t$ 
points in $\R^3$:
\begin{eqnarray}
C_{n,t}(\phi) & = &
\{ (s_1, \ldots , s_n, x_1, \ldots , x_t) 
\in {(S^1)}^{n}\times (\R^3)^t  |  \nonumber\\
& & s_1 < \cdots < s_n ,\ {\mathrm cyclically\ ordered\ on\ the\ circle}, 
\label{Cnt}\\
& &  x_i\neq x_j \ {\mathrm if}\ i\neq j, \nonumber\\
& &  {\phi}(s_i) \neq x_j \ {\mathrm if} \ i\in \{1,\ldots , n\},\,
j\in \{1,\ldots , t\} \} . \nonumber
\end{eqnarray}
Let us pick two bijections $\beta : V_\Gamma^i \rightarrow \{1,\ldots,t\}$,
and $\gamma : V_\Gamma^o \rightarrow \{1,\ldots,n\}$, such that 
if $\gamma(u)=i$, and $v$ is the next vertex encountered 
when going around $W_\Gamma$ according to its orientation,
then $\gamma(v)=i+1 \bmod n$.
Every point of $C_{n,t}(\phi)$ corresponds to an embedding
$\Phi : \Gamma \rightarrow\R^3$ (see fig. \ref{huit}), defined by 
\begin{equation}
\Phi(v) = \left\{
\begin{array}{ll}
x_{\beta(v)} & \ {\mathrm if}\ v\in V_\Gamma^i, \\
\phi(s_{\gamma(v)}) & \ {\mathrm if}\ v\in V_\Gamma^o,
\end{array} \right.
\label{plonge}
\end{equation}
such that $\Phi(e)=\phi([s_i,s_{i+1}])$ if $e\subset W_\Gamma$,
$\partial e = (u,v)$, $\gamma(u)=i$, $\gamma(v)=i+1$, and
if $e\in E_\Gamma^{nW}$, $\Phi(e)$ is the straight line segment
joining $\Phi(u)$ to $\Phi(v)$, with $\partial e = (u,v)$.
Hence we can identify $W_\Gamma$ and $S^1$ via $\gamma$, and
think of $\Phi$ as an extension of $\phi$ to $\Gamma$.
Clearly, $\Phi$ depends on the choice of $\beta,\gamma$, however
$I(\Gamma)$ doesn't.

In (\ref{defomega}) we omitted the maps $\beta,\gamma$ from the notations,
thus the indices $i,j$ in $x_i,s_j$ correspond to certain
specific vertices in $V_\Gamma^i,V_\Gamma^o$.

Now, if we choose an orientation $\orga$ for the edges of 
$\Gamma$, we can define a continuous form on $C_\Gamma$ 
of degree $2 |E_{\Gamma}^{nW}|$  ($|E_{\Gamma}^{nW}|$ 
is the number of internal edges)
as follows:
\begin{equation}
\omega(\Gamma,\orga) =\prod_{e\in E_{\Gamma}^{nW}} \phi^*_e \omega,
\end{equation}
where $\omega$ is the Gauss two-form on $S^2$:
\begin{equation}
\omega(x) = {1\over 4\pi} (x^1 dx^2 \wedge dx^3 +
             x^3 dx^1 \wedge dx^2 + x^2 dx^3 \wedge dx^1),
\end{equation}
$\phi_e : C_\Gamma \rightarrow S^2$
is given by
\begin{equation}
\left\{
\begin{array}{ll}
u(x_j-x_i) & \mbox{if $e$ is an oriented 
edge  $\partial e=(i,j)$ } \\
&\mbox{connecting two distinct internal vertices $i,j \in V_\Gamma^i$ },\\
u(\phi(s_j)-x_i)
&\mbox{if $e$ is an oriented edge $\partial e=(i,j)$}\\
&\mbox{connecting an internal vertex $i$ with an external vertex $j$,}\\
u(\phi(s_j)-\phi(s_i))
&\mbox{if $e$ is an oriented edge $\partial e=(i,j)$ }\\
&\mbox{connecting two distinct external vertices $i,j \in V_\Gamma^i$,}
\end{array}
\right.
\label{defomega}
\end{equation}
and $u : \R^3-\{0\} \rightarrow S^2$ is defined by 
\begin{equation}
u(x)=\frac{x}{|x|}
\label{defu}
\end{equation}
for $x=(x^1,x^2,x^3)\in\R^3$.

\noindent
{\bf Remarks:}
\begin{enumerate}
\item We gave the value of $\phi_e$ only in the interior of 
$C_\Gamma$. The key point is that the fiber of $C_\Gamma$,
which is a compact manifold with corners, 
is the compactification of $C_{n,t}(\phi)$, on which
$\phi_e$ admits a smooth extension.
\item $\omega(\Gamma,-\orga) = -\omega(\Gamma,\orga) =
- \omega(-\Gamma,\orga)$.
\item If $\Gamma$ contains an edge $e\in E^{nW}_\Gamma$
whose boundary consists of only one vertex, 
$\partial e=(i,i)$, $i\in V_\Gamma$, a case which is not covered by
(\ref{defomega}), it will be convenient to set $\omega(\Gamma,\orga)=0$.
\end{enumerate}
Then if we choose an orientation $\Omega$ on the fiber of 
$C_\Gamma$, we can set
\begin{equation}
I(\Gamma,\orga,\Omega) =(-1)^{{\mathrm deg}(p_\Gamma)_* \omega}\, 
(p_\Gamma)_*  \omega(\Gamma,\orga),
\end{equation}
where $(p_\Gamma)_*$ is the push-forward (integration along the fiber),
whose definition is given in \ref{appa}.
(An orientation is needed in order to integrate forms.)
With this definition we see that 
\begin{equation}
I(\Gamma,\orga,\Omega) =I(-\Gamma,\orga,\Omega) =
-I(\Gamma,-\orga,\Omega) =-I(\Gamma,\orga,-\Omega).
\end{equation}
The last two equalities are obvious, the first one follows from the fact 
that $\omega(\Gamma,\orga)$ does not depend on the vertex orientation.
\begin{prop}
If $\Gamma$ is a Wilson V-oriented graph,
which is either trivalent, or trivalent except for one tetravalent vertex,
and  $\orga$ is an orientation of its edges, then there exists
an orientation $\Omega(\Gamma,\orga)$ on the fiber of $C_\Gamma$ such that 
\begin{equation}
\Omega(\Gamma,\orga)=-\Omega(\Gamma,-\orga)
= -\Omega(-\Gamma,\orga).
\label{omor}
\end{equation}
\label{promor}
\end{prop}
Assuming the proposition, we put
\begin{equation}
I(\Gamma)=I(\Gamma,\orga,\Omega(\Gamma,\orga)),
\end{equation}
and then (\ref{reverse}) follows immediately. 

It is possible to prove (\ref{omor}), and thus to
define $I(\Gamma)$, without any assumption on the
valences, but we shall refrain from doing that since we don't need
to deal with with arbitrary graphs in this paper. This would lead
to the full definition of the graph complex.

In the case of trivalent graphs, $I(\Gamma)$ is a function on 
$\cal K$, and if $\phi \in \cal K  $
\begin{equation}
I(\Gamma)(\phi) =  \int_{C_{n,t} (\phi) } \omega (\Gamma, \orga),
\end{equation}
where the orientation on $ C_\Gamma (\phi) $ 
is given by $\Omega(\Gamma,\orga)$.
The existence of a compactification of 
$C_\Gamma (\phi)$ on which $ \omega (\Gamma, \orga) $
extends smoothly tells us that this integral is convergent. 

{\em Proof of proposition \ref{promor}.}
Let $\phi \in \cal K$ be a circle embedding,
and $\Gamma$ a trivalent Wilson V-oriented graph equipped with an 
orientation of the edges $\orga$. 
The interior $C_{n,t}(\phi)$ of the fiber of $C_\Gamma$ over $\phi$
(see (\ref{Cnt})) is included in
$ {S_1}^{V_{\Gamma}^{o}} 
\times ({\R^3})^{V_{\Gamma}^{i}}$. Here ${S_1}^{V_{\Gamma}^{o}}$ is  
the set of maps from ${V_{\Gamma}^{o}}$ to $S_1$.
Coordinates on this space are denoted by $ X_v^i$, $v\in V_\Gamma$, 
$i\in\{1,2,3\}$. If $v\in V_\Gamma^o$, $i=1$ and $X_v^1\in S^1$.
If $v\in V_\Gamma^i$, $X_v=(X_v^1,X_v^2,X_v^3)\in\R^3$.
Remember the maps $o_v, v\in V_\Gamma$ defining the 
V-orientation (see \ref{ov}): $o_v =1 \mbox{ if } v \in V_\Gamma^o $,
$o_v \in \{1,2,3\} \mbox{ if } v \in V_\Gamma^i $ is trivalent.
In the particular case of {\it trivalent} graphs 
we define a volume form as follows:
\begin{equation}
\Omega(\Gamma,\orga) = 
\bigwedge_{ {e \in E_{\Gamma}^{nW}}}\Omega_e,
\end{equation} 
\begin{equation}
\Omega_e = dX_v^{o_v(e)} \wedge dX_u^{o_u(e)},
 \mbox{ if } {\partial e =(u,v)}.
\end{equation}
It is straightforward to check that $\Omega(\Gamma,\orga)$ 
does not depend on the choice of maps $o_v,v\in V_\Gamma$,
but only on their oriented class.

Now let $\Gamma_x$ be a V-oriented graph with one tetravalent vertex $x$.
Suppose first that $x$ is an internal vertex, 
and 
let $e_1,e_2,e_3,e_4$ be the four edges surrounding $x$, with 
$\partial e_i =(x,v_i)$.
The labels are chosen such that $ o_x(e_i)=i$.
Define $\Omega(\Gamma_x,or_{\Gamma_x}) =\Omega_1 \wedge \Omega_2 $
with
\begin{equation}
\Omega_1 = d^3X_x \wedge dX_{v_1}^{o_{v_1}(e_1)}
\wedge dX_{v_2}^{o_{v_2}(e_2)} \wedge dX_{v_3}^{o_{v_3}(e_3)} 
\wedge dX_{v_4}^{o_{v_4}(e_4)},
\end{equation}
where $d^3 X =dX^1 \wedge dX^2 \wedge dX^3$, and 
\begin{equation}
\Omega_2 = \bigwedge_{e \in E^{nW}_{\Gamma_x} \atop 
\neq {e_1,e_2,e_3,e_4}} \Omega_e.
\end{equation}
Next if $x$ is a vertex on the Wilson line,
let $e_1,e_2$ be the two internal edges with 
$\partial e_i =(x,v_i)$ and 
$w_1,w_2$ the two Wilson edges meeting the vertex $x$.
The labels are chosen such that : $ o_x(e_i)=i$ and 
$ \partial w_1 =( \cdot , x), \partial w_2 =(  x, \cdot)$.
Define
$\Omega(\Gamma_x,or_{\Gamma_x}) =\Omega_1 \wedge \Omega_2 $
with
\begin{equation}
\Omega_1 = dX_x \wedge dX_{v_1}^{o_{v_1}(e_1)}
\wedge dX_{v_2}^{o_{v_2}(e_2)} ,
\end{equation}
 and 
\begin{equation}
\Omega_2 = \bigwedge_{e \in E^{nW}_{\Gamma_x} \atop 
\neq {e_1,e_2}} \Omega_e.
\end{equation}
By construction  $\Omega(\Gamma_x,{{\mathrm or}\,}_{\Gamma_x})$ 
does not depend on the choice of maps $o_v,v\in V_\Gamma$,
but only on their oriented class. 

With these definitions, the verification of (\ref{omor}) presents
no difficulty.
$\Box$

\section{Invariance under $\delta$}
In this section we denote by $|\Gamma|$ and $|\Gamma_+|$ the
orders of the groups ${\mathrm Aut}(\Gamma)$ and ${\mathrm Aut}_+(\Gamma)$.
Let $G^3_n =\{ \Gamma \in G^3 | \mbox{deg }\Gamma = n\}$, and 
\begin{equation}
Z_n = \sum_{\Gamma\in G^3_n }{1 \over |\Gamma|}D(\Gamma) I(\Gamma).
\end{equation}
Thus $Z_n$ is a function on $\cal{K}$ with values in $ {\cal{A}}_n$
(the space of BN diagrams of order $n$).
We remark that it is independent, 
as the combination $D(\Gamma) I(\Gamma)$ is, of the choice 
of vertex orientation for each trivalent Wilson graph.
Another observation is that $I(\Gamma)/|\Gamma|$ is a sum over all
the embeddings (\ref{plonge}) of $\Gamma\in G^3$ in $\R^3$, weighted by factors
$\omega(\Gamma,\orga)$.

\begin{figure}[htbp]
$$\epsf{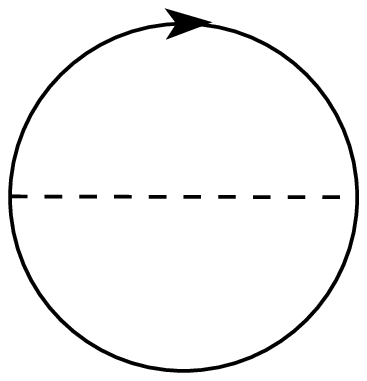}[xscale=1/2,yscale=1/2] 
$$ \caption{\label{theta}}
\end{figure}

\noindent
The expectation value of a Wilson loop is a formal power series in 
$\hbar$ with values in $\cal{A}$:
\begin{equation}
Z = 1 + \sum_{n=1}^\infty \hbar^n Z_n.
\end{equation}
We define 
\begin{equation}
\delta I(\Gamma)  = \sum_{e \in E^a_\Gamma }  I(\delta_e\Gamma),
\end{equation}
where the sum is over all admissible edges of $\Gamma$. 
Observe that $\delta I(\Theta)=0$ according to our Feynman rules,
where $\Theta $ is the unique Wilson trivalent graph of order $1$ 
represented in fig. \ref{theta}.
Then we have the following theorem, which has been known for
some time, see e.g. \cite{kohno} for a discussion
in the setting of the whole graph complex.
\begin{theo}
\begin{equation}
\delta Z = 0.
\end{equation}
\label{delta}
\end{theo}
{\em Proof.}
We define the following equivalence relation on $E^a_\Gamma$ :
\begin{equation}
 e \sim f \mbox{ iff }  \delta_e \Gamma \sim 
 \pm \delta_f \Gamma
\end{equation}
and denote by ${\cal E}_\Gamma = E^a_\Gamma/ \sim$  the set of equivalence
classes.
Let $ S^\pm_{\Gamma,e} =\{ f\in E^a_\Gamma | \delta_f \Gamma \sim 
 \pm \delta_e \Gamma \} $.
Observe that if $S^+_{\Gamma,e} \bigcap S^-_{\Gamma,e} \neq \emptyset $,
then $I(\delta_e\Gamma)=0$.
Then using the definition of $\delta$ and the latter property 
we clearly have:
\begin{equation}
\delta Z_n= 
\sum_{\Gamma\in G^3_n } { 1 \over |\Gamma|} D(\Gamma) 
\sum_{e\in {\cal E}_\Gamma} I(\delta_e\Gamma)
(|S^+_{\Gamma,e}| -|S^-_{\Gamma,e}|)
\end{equation}

\begin{figure}[htbp]
$$\epsf{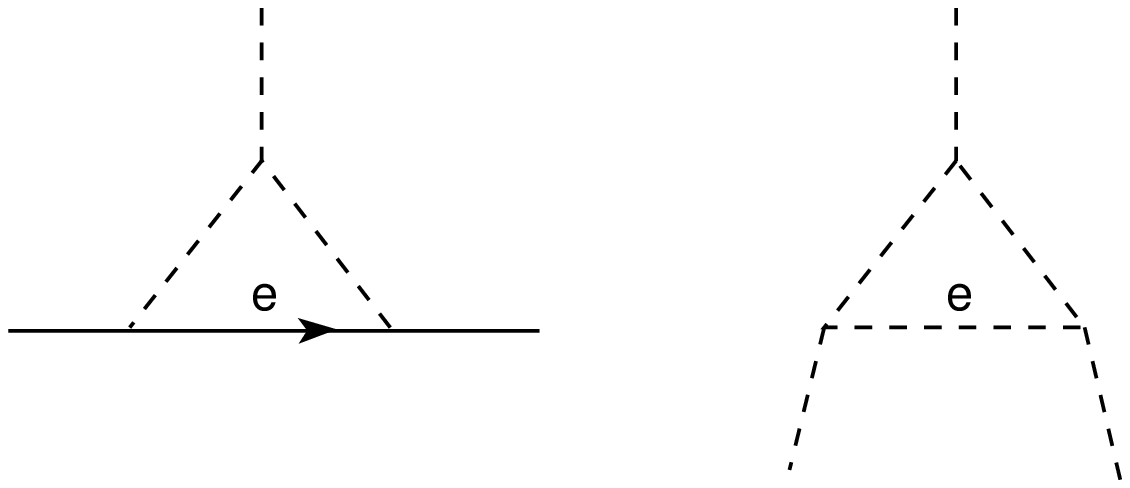}[xscale=2/3,yscale=2/3] 
$$ \caption{\label{triangle}}
\end{figure}

\begin{figure}[htbp]
$$\epsf{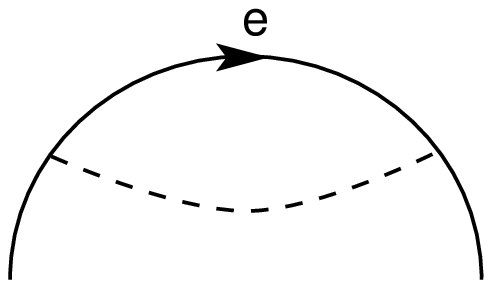}[xscale=2/3,yscale=2/3] 
$$ \caption{\label{triangle2}}
\end{figure}

\begin{lemma}
If $e\in E^a_\Gamma$ and if e is not as in figs.
\ref{triangle},\ref{triangle2},
then 
\begin{equation}
|S^\pm_{\Gamma,e}| = m_\pm (\Gamma,e) 
{ |\Gamma_+|\over |\delta_e \Gamma_+ |}
\end{equation}
where $m_\pm (\Gamma,e) = \# \, i \in \{1,2,3\}$ such that
$\Gamma_x^{(i)} \sim \pm \Gamma$,
with $\Gamma_x = \delta_e \Gamma$.
\label{S}
\end{lemma}

Assuming this lemma, the proof of theorem \ref{delta} goes as follows.
First we can assume that $ |\Gamma_+|=|\Gamma|$  and 
$|\delta_e\Gamma_+| =|\delta_e \Gamma| $, otherwise 
$D(\Gamma) =0$ or $ I(\delta_e \Gamma) =0$.
Let us denote by $ {\cal E}^\prime_\Gamma$
the subset of 
${\cal E}_\Gamma $ 
which satisfies the assumption of the lemma.
Then 
\begin{equation}
\delta Z_n =
\sum _{\Gamma \in G_n^3}  D(\Gamma) \sum_{e\in {\cal E}^\prime_\Gamma }
I(\delta_e \Gamma) {1\over |\delta_e\Gamma|} ( m_+(\Gamma,e)- m_-(\Gamma,e)),
\end{equation}
because $ I(\delta_e \Gamma) =0 $ if $e$ is as in figs. \ref{triangle}
or \ref{triangle2}.
By lemma \ref{croix} the coefficient of $I(\Gamma_x)$, where 
$\Gamma_x =\delta_e \Gamma$, in $\delta Z_n$
is 
\begin{equation}
C(\Gamma_x) = {1\over |\Gamma_x| }
\sum D(\Gamma_x^{(i)}) ( ( m_+(\Gamma_x^{(i)},e)- m_-(\Gamma_x^{(i)},e)),
\end{equation}
where the sum is over the inequivalent $\Gamma_x^{(i)}$ such that
$\delta_e\Gamma_x^{(i)}=\Gamma_x$.
The relations :
\begin{equation}
\begin{array}{c}
 D(\Gamma_x^{(1 )})+D(\Gamma_x^{( 2)})+D(\Gamma_x^{( 3)}) = 0 \\
 D(\Gamma) + D(-\Gamma) = 0
\end{array}
\label{stu}
\end{equation}
imply that $C(\Gamma_x)=0$, which concludes the proof of theorem 
\ref{delta}. 
This can be shown by considering the three cases :
 \begin{enumerate}
\item All $\Gamma_x^{(i)},\, i= 1,2,3$ are inequivalent. \\
Then  $ m_+(\Gamma_x^{(i)},e) =1 \,\forall i$, and  
$ m_-(\Gamma_x^{(i)},e) =0 \,\forall i$.
Thus $C(\Gamma_x)=0$ by the IHX (STU) relations (\ref{stu}). 
 \item  $\Gamma^{(1)} \sim \pm \Gamma^{(2)} \neq \Gamma^{(3)} $. \\
Then $ m_+(\Gamma_x^{(3)},e) =1$,  
$ m_-(\Gamma_x^{(3)},e) =0 $ \\
and $ m_+(\Gamma_x^{(1)},e) \in \{1,2\} $, 
$ m_-(\Gamma_x^{(1)},e) \in \{0,1,2\}$. \\
We can assume $  m_+(\Gamma_x^{(1)},e) + m_-(\Gamma_x^{(1)},e) = 2$, 
otherwise $\Gamma^{(1)} \sim - \Gamma^{(1)}$, so \\
$ D(\Gamma^{(1)}_x) =D(\Gamma^{(2)}_x)=0$, 
and $C(\Gamma_x)=D(\Gamma^{(3)}_x)=0$
by IHX. \\
If $ m_-(\Gamma_x^{(1)},e) =0,\,  m_+(\Gamma_x^{(1)},e) =2 $,
then $  D(\Gamma^{(1)}) =D(\Gamma^{(2)})$ and \\
$  C(\Gamma_x)= 2D(\Gamma^{(1)}) +D(\Gamma^{(3)})=0$
by IHX. \\
If $ m_-(\Gamma_x^{(1)},e) =1,\,  m_+(\Gamma_x^{(1)},e) =1 $,
then $  D(\Gamma^{(1)}) +D(\Gamma^{(2)})=0$ and 
$  C(\Gamma_x)= D(\Gamma^{(3)})=0$
by IHX.
\item  $\Gamma^{(1)} \sim \pm \Gamma^{(2)} \sim \pm \Gamma^{(3)} $.
Then  $ D(\Gamma^{(i)})=0\, \; , i= 1,2,3 $,
and $C(\Gamma_x) =0$ because it is proportional to  $D(\Gamma^{(i)})$. 
$\Box$
\end{enumerate}
{\em Proof of lemma \ref{S}.}
Let 
\begin{equation}
\begin{array}{ll}
O_g^\pm (\Gamma,e) &= \{(\rho,\sigma)\in {\mathrm Aut}(V_\Gamma) 
\times {\mathrm Aut}(E_\Gamma)
 |( {\rho},{\sigma}) \cdot \delta_g \Gamma =\pm \delta_e \Gamma \}, \\
P_i^\pm (\Gamma,e) &= \{ 
({ \rho},{ \sigma})\in {\mathrm Aut}(V_\Gamma) \times {\mathrm Aut}(E_\Gamma)
|({ \rho},{ \sigma}) \cdot \Gamma =\pm \Gamma_x^{(i)}, \;\;
\Gamma_x=\delta_e \Gamma \},
\end{array}
\end{equation}
where $g\in E_\Gamma^a$, $i\in \{1,2,3\}$.
By proposition 1
\begin{equation}
|\bigcup_{g \in E^a_\Gamma} O_g^\pm (\Gamma,e)| 
=| \bigcup_{i \in \{1,2,3\} } P_i^\pm (\Gamma,e)|.
\label{UU}
\end{equation}
There is an obvious action of ${\mathrm Aut}_+ (\delta_e\Gamma)$ on 
$O_g^\pm (\Gamma,e)$ which is free and transitive.
Therefore, either $O_g^\pm (\Gamma,e)= \emptyset $
or $O_g^\pm (\Gamma,e) \sim {\mathrm Aut}_+ (\delta_e\Gamma)$ as a set.
By lemma \ref{lemef}, we have 
\begin{equation}
g \in S^\pm_{\Gamma,e} \Leftrightarrow O_g^\pm (\Gamma,e)\neq \emptyset, 
\end{equation}
so $|S^\pm_{\Gamma,e}| = \# g \in E_\Gamma^a \mbox{ s.t. } 
O_g^\pm(\Gamma,e)\neq \emptyset $,
and since $ O_g^\pm (\Gamma,e) \bigcap O_h^\pm (\Gamma,e) =\emptyset $
if $g\neq h$, we get:
\begin{equation}
|\bigcup_g O_g^\pm (\Gamma,e)| = |S^\pm_{\Gamma,e}||\delta_e\Gamma_+|.
\label{O}
\end{equation}
Similarly, ${\mathrm Aut}_+(\Gamma)$ acts freely and transitively on 
$P_i^\pm (\Gamma,e)$ for each $i \in \{1,2,3\}$. 
Moreover we have the following property, whose proof is left to the reader:
\begin{lemma}
If $e$ satisfies the assumption of lemma \ref{S}, and $i\neq j$,
\begin{equation}
 P_i^\pm (\Gamma,e) \bigcup P_j^\pm (\Gamma,e) = \emptyset. 
\end{equation}
\end{lemma}
This implies 
\begin{equation}
| \bigcup_i P_i^\pm (\Gamma ,e)| =m_\pm (\Gamma ,e) |\Gamma_+|.
\label{P}
\end{equation}
The equalities (\ref{UU}),(\ref{O}), and (\ref{P}) imply the lemma \ref{S}. 
$\Box$

\section{Anomalies}
Let $\Gamma$ be a V-oriented Wilson trivalent graph,
and denote by $d$ the exterior differential on $\cal K$.
Using the definition of $I(\Gamma)$ as a pushforward along the 
compact fiber of $C_\Gamma$, and the commutation relation 
between $d$ and the pushforward given in (\ref{diff}),
the variation of $I(\Gamma)$ under a change of embedding
can be expressed as a sum over all the 
strata  of  $C_\Gamma$ \cite{bt}:
\begin{equation}
dI(\Gamma) = \sum_{e \in E^a_\Gamma } 
\epsilon(\delta_e{\Gamma}) I(\delta_e\Gamma)
 +\delta_a I(\Gamma),
\end{equation}
where $\epsilon(\delta_e{\Gamma})I(\delta_e\Gamma)$ 
corresponds to the pushforward along the codimension one strata 
$\partial_e C_\Gamma $, when $e$ is not as in fig. \ref{triangle2},
and is obtained by collapsing two vertices of $\Gamma$
along the admissible edge $e$. 
The signs $\epsilon(\delta_e{\Gamma})=\pm 1$ 
depend on the induced orientations of these strata.
The ``anomalous'' term $\delta_a \Gamma$ is the contribution of
all the other strata.
Bott and Taubes showed that the pushforward along 
these strata is in general zero except for special ones, which 
correspond, in the case of prime graphs, to the simultaneous 
collapse of all the vertices together.
Moreover, there exists a ``universal'' way of calculating the contribution
of these strata (universal means here independent of the embedding).
More precisely, we can state their result as follows:
if $\Gamma$ is a {\em prime} V-oriented Wilson trivalent graph,
then
\begin{equation}
\delta_a I(\Gamma) =  
f_{\Gamma} \, { dI(\Theta)  \over 2}.
\label{ano}
\end{equation}
Here $dI(\Theta)$ is the differential of the self-linking integral $I(\Theta)$:
if $\phi $ is an embedding of $S^1$ into $\R^3$ then 
\begin{equation}
I(\Theta) (\phi) ={1\over 4\pi} \int_{S^1 \times S^1}  ds_1 ds_2
\,\epsilon_{\mu \nu \rho}
 \dot{\phi}^\mu(s_1)\dot{\phi}^\nu(s_2) 
{ \phi^\rho(s_2)-\phi^\rho(s_1) \over |\phi(s_2)-\phi(s_1)|^3}.
\end{equation}
The constant of proportionality $f_{\Gamma}$ is independent of the 
embedding $\phi$ and is expressed as an integral:
\begin{equation}
f_\Gamma =\int_{S{n,t}}\Theta_{\Gamma},
\end{equation}
where $S_{n,t}$ is the following variety of dimension $n+3t$,
$n$ being the number of Wilson vertices and $t$ the number 
of internal vertices: it is the set of
$(a,\eta_1, \ldots , \eta_n, \omega_1, \ldots , \omega_t) 
\in S^2 \times {\R}^{n}\times ({\R}^3)^t$
such that:
\begin{equation}
\label{Snt}
\begin{array}{l}
\eta_1 < \cdots < \eta_n , \mbox{ or cyclic permutations}\\
\omega_i\neq \omega_j \mbox{ if } i\neq j, \\
a\cdot \eta_i \neq \omega_j \mbox{ if } i\in \{1,\ldots , n\},\,
j\in \{1,\ldots , t\}, \\
\sum_{i=1}^n \eta_{i}^2 + \sum_{i=1}^t |\omega_{i}|^2 =1,\\
\sum_{i=1}^n  \eta_{i} + \sum_{i=1}^t {<\omega_{i},a>} =0.\\
\end{array}
\end{equation}
If $\Gamma$ is equipped with an orientation of its edges and an 
orientation of its vertices, one defines an orientation 
$\Omega_{n,t}$ of $S_{n,t}$ 
induced by $\Omega(\Gamma,or_\Gamma)$,
and a $n+3t$-form 
$\Theta_{\Gamma}$:
\begin{equation}  
\Theta_{\Gamma}=\prod_{e\in E_{\Gamma}^{nW}} \Theta_e,
\end{equation}
the product being over all the internal oriented edges of $\Gamma$, and 
$\Theta_e = \phi^*_e \omega$,
where $\omega$ is the Gauss two-form, and
$\phi_e : S_{n,t} \rightarrow S^2$ is given by
\begin{equation}
\left\{
\begin{array}{ll}
u(\omega_j-\omega_i) & \mbox{if $e$ is an oriented 
edge  $(i,j)$ connecting two internal vertices $i$ and $j$,}\\
u(a\cdot \eta_j-\omega_i) 
&\mbox{if $e$ is an oriented edge $(i,j)$ connecting an internal vertex $i$}
\\ & \mbox{with an external vertex $j$,}
\\
u(a(\eta_j-\eta_i))
&\mbox{if $e$ is an oriented edge $(i,j)$ connecting 
two external vertices $i$ and $j$,}
\end{array}
\right.
\end{equation}
where $u$ is defined by (\ref{defu}).

Note that a simple computation gives $f_\Theta =2$, which is consistent
with our definition of $\delta$ and (\ref{diffe}).

A priori the integral $f_\Gamma$ is not well-defined because of the
singularities at coinciding points. We should verify that the integral
is indeed convergent. This follows from the fact that $S_{n,t}$ admits
a compactification, which is a manifold with corners, such that
$\Theta_{\Gamma}$ extends smoothly on it.

We would like to emphasize that 
the definitions of the V-orientation and $\Omega(\Gamma,\orga)$
we gave imply that $\epsilon(\delta_e{\Gamma})=+1$, hence
\begin{prop}
\label{proepsi}
The total variation of $I(\Gamma)$ is:
\begin{equation}
dI(\Gamma) = \sum_{e \in E^a_\Gamma }   I(\delta_e\Gamma)
 +\delta_a I(\Gamma).
\label{diffe}
\end{equation}
\end{prop}
{\em Proof.}
Let $\Gamma$ be a trivalent graph, $e \in E_\Gamma^a$ an admissible edge,
$\bar{C}_{n,t}(\phi)$ the fiber of $C_\Gamma$ over 
$\phi \in \cal K$, and $\partial_e \bar{C}_{n,t}(\phi)$ the codimension 
one strata of $\bar{C}_{n,t}(\phi)$ corresponding to the 
collapse of $e$.
Coordinates in $ C_{n,t}(\phi)$ can be taken 
to be 
\begin{equation}
((s_i)_{i=1,\cdots,n }, (X_i)_{i=1,\cdots,t})
\in   (S_1)^{n} \times (\R^3)^t
\end{equation}
If $e$ is an edge on the Wilson line $\partial e =(i+1,i)$,
the strata associated to $e$ corresponds to $s_i=s_{i+1}$, and 
the in-going normal vector field to this stratum is
\begin{equation} 
n_e = {\partial \over \partial s_{i+1} } -{\partial \over \partial s_i},
\end{equation}
since $s_{i+1} > s_i$ in $C_{n,t}(\phi)$.
It is clear that 
\begin{equation}
i_{n_e} \Omega(\Gamma,\orga) = 
\Omega(\delta_e\Gamma,or_{\delta_e\Gamma}).
\end{equation}
If $e$ is an internal edge connecting two internal vertices 
 $\partial e =(i,j)$,
the vicinity of the stratum associated to $e$ 
corresponds to 
$X_j =X + {\alpha \over 2}u,\,X_i =X - {\alpha \over 2}u $,
where $\alpha>0$ and $u \in S^2$.
So the in-going normal vector field to this stratum is
$n_e = {\partial \over \partial \alpha }$ 
and  
$i_{n_e} \Omega(\Gamma,\orga) = \alpha^2 \omega(u) \wedge
\Omega(\delta_e\Gamma,{{\mathrm or}\,}_{\delta_e\Gamma})$
($\omega $ is the Gauss form).
Then using the fact that 
$\partial_e C_{n,t}(\phi) = S^2 \times  C_{n,t-1}(\phi)$
and integrating over $S^2$, 
the induced orientation on $C_{n,t-1}(\phi)$
is $\Omega(\delta_e\Gamma,or_{\delta_e\Gamma})$.

If $e$ is an internal edge connecting one internal vertex
and one Wilson vertex a similar analysis can be performed. 
$\Box$

\noindent
The vanishing theorems of Bott and Taubes can be improved 
by showing that:
\begin{theo}
\label{annul}

\begin{description}
\item[(i)] If $\Gamma$ is not primitive then $f_\Gamma =0 $. 
\item[(ii)] If $\Gamma$ is primitive and ${\mathrm deg}(\Gamma)$ is even 
then $f_\Gamma =0$.
\end{description}
\end{theo}
{\em Proof of (i).} 
Denote by $\Gamma_\alpha $ the connected subgraphs of $\Gamma - W_\Gamma$, 
\begin{equation}
\Gamma - W_\Gamma =\bigcup_{\alpha=1}^n \Gamma_\alpha.
\end{equation}
Consider the following vector fields on 
$\R^n \times (\R^3 )^t$: 
\begin{equation}
V_\alpha = 
\sum_{i\in V_\Gamma^o\cap\Gamma_\alpha}
{\partial \over \partial \eta_i}
+
\sum_{j \in V_\Gamma^i\cap\Gamma_\alpha} 
\langle a \, , \, {\partial \over \partial \omega_j}\rangle,
\end{equation}
\begin{equation}
{\cal D} = 
\sum_{i\in V_\Gamma^o} 
\eta_i{\partial \over \partial \eta_i}
+
\sum_{j \in V_\Gamma^i}
\langle \omega_j \, , \, {\partial \over \partial \omega_j}\rangle.
\end{equation}
$V_\alpha $ is the vector field of translation of all 
vertices of one subgraph $\Gamma_\alpha$, and $\cal D$ is the vector
field of global dilatation. Now put
\begin{equation}
V=
\sum_\alpha c_\alpha V_\alpha  -\lambda{\cal D},
\end{equation}
where $c_\alpha$ are not all simultaneously zero, and satisfy the condition
\begin{equation}
\sum_\alpha c_\alpha {\mathrm deg}\,\Gamma_\alpha =0,
\end{equation}
and $\lambda$ is the following function on $S_{n,t}$ :
\begin{equation}
\lambda = \sum_\alpha c_\alpha 
\left ( \sum_{i\in V_\Gamma^o \cap \Gamma_\alpha} \eta_i
+ \sum_{j \in V_\Gamma^i \cap \Gamma_\alpha} 
\langle a \, , \, \omega_j\rangle\right).
\end{equation}
Then part (i) of the theorem 
is a direct consequence of the following two properties:
\begin{equation}
\begin{array}{c}
\mbox{If $\Gamma$ is not primitive, } V \in TS_{n,t}, \, V\neq 0 
\mbox{ almost everywhere.} \\
i_V \Theta_\Gamma =0.
\end{array}
\end{equation}
Here $i_V$ denotes the interior product.
The first property follows from the fact that $V$ preserves the 
conditions defining the embedding
of $S_{n,t}$ into $S^2 \times {\R}^{n}\times ({\R}^3)^t$, 
the second from ${\cal D}(\phi_e)=V_\alpha (\phi_e)=0$. $\Box$ 

\noindent
{\em Proof of (ii).}
Consider the diffeomorphism $S$
of $S_{n,t}$ defined by:
\begin{equation}
\begin{array}{c}
S(a)=-a,\\
S(\eta_i)=\eta_i,\\
S(\omega_j)= -\omega_j.
\end{array}
\end{equation}
We find that
the behaviour of $\Theta_{\Gamma}$ and the orientation of 
$S_{n,t}$ is given by:
\begin{equation}
\begin{array}{c}
S^{*}\Theta_{\Gamma}= (-1)^{n+3t\over 2}\Theta_{\Gamma},\\
S^{*}\Omega_{n,t} =  (-1)^{t+1}\Omega_{n,t}. 
\end{array}
\end{equation}
Thus 
\begin{equation}
\int_{S{n,t}} \Theta_{\Gamma} = (-1)^{n+3t\over 2}\int_{S{n,t}}
 S^*\Theta_{\Gamma}
=(-1)^{{n+t\over 2} +1}\int_{S{n,t}} \Theta_{\Gamma},
\end{equation}
and $f_{\Gamma}=0$ if the graph is of even order
($(n+t)/2$ is the order of the graph). $\Box$

\begin{figure}[htbp]
$$\epsf{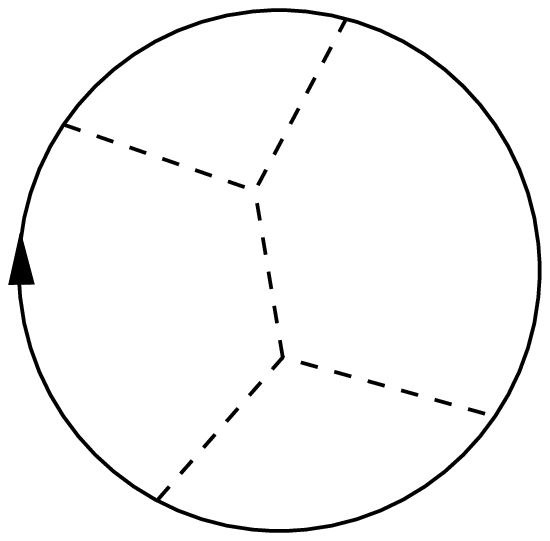}[xscale=2/3,yscale=2/3] 
$$ \caption{\label{YY}}
\end{figure}

So the first non-trivial prime graph which could be anomalous,
apart from $\Theta$, appears at order three 
(this is the only one at this order) and is given in fig. \ref{YY}. 

Combining (\ref{diffe}) with theorem \ref{delta},
we can now state the invariance
theorem of the expectation value of a Wilson loop.
First, let us denote by $\alpha$ the element of $\cal A$ 
characterizing the anomalous variations :
\begin{equation}
\alpha = {1 \over 2} \sum_{\Gamma\;{\mathrm primitive}}
 \hbar^{{\mathrm deg}\Gamma}
f_\Gamma 
{D(\Gamma) \over |\Gamma|}.
\end{equation}

The self-linking integral $I(\Theta) (\phi)$ is not an invariant,
but it is well-known \cite{calu} that if we introduce a framing, given by a 
normal vector field $\nu$, then the linking number of the two curves
$\phi$ and $\phi+\nu$ is
\begin{equation}
{\mathrm lk}(\phi,\phi+\nu) = I(\Theta) (\phi) + \tau(\phi,\nu),
\label{tors}
\end{equation}
where $\tau(\phi,\nu)$ is the total torsion, given by
\begin{equation}
\tau(\phi,\nu) = {1\over 2\pi}\int ds\, \dot{\phi}(s) { (\dot{\phi}(s),\nu(s),
\dot{\nu}(s)) 
\over |\dot{\phi}(s)\wedge \nu(s)|^2 }.
\end{equation}
This implies
\begin{theo}
\label{invariant}
\begin{equation}
\hz(\phi,\nu) = Z(\phi) \exp
\left( \alpha \, \tau(\phi,\nu) \right)
\end{equation}
is a framed knot invariant.
\end{theo}
In the next section, we are going to show that it is in fact a 
universal Vassiliev invariant.
The theorem is a consequence of the following important property of $Z$:
\begin{lemma}
Let $G^3_*$ denote the trivalent graphs $\Gamma$ with 
${\mathrm deg}\Gamma > 0$. Then
\begin{equation}
\log Z =  \sum_{\Gamma \in G^3_*} \hbar^{{\mathrm deg}\Gamma} 
{I({\Gamma}) \over {| \Gamma |}} C(D(\Gamma)).
\end{equation}
\label{grlike}
\end{lemma}
By lemma \ref{lemC}, the map $C$ appearing in this
equality is a projector on the
subspace of primitive elements of $\cal A$.
Hence $\log Z$ is a primitive element of
the completion $\hat{\cal A}$ of the Hopf algebra $\cal A$. 
Thus\footnote{We owe this remark to T. T. Q. Le.}
lemma \ref{grlike} also shows:
\begin{theo}
$\hz(\phi,\nu)$ is a group-like element of $\hat{\cal A}$.
\end{theo}
It is known \cite{bng,tle} that the universal Vassiliev invariant
constructed from the KZ connection satisfies the same group-like property. 

Assuming lemma \ref{grlike}, the proof of theorem \ref{invariant}
is as follows:
by theorem \ref{delta}, part(i) of theorem \ref{annul} and (\ref{ano}),
\begin{equation}
d \log Z = {d I(\Theta)\over 2} 
\sum_{\Gamma\;{\mathrm primitive}}
\hbar^{{\mathrm deg}\Gamma} f_{\Gamma} {D(\Gamma) \over {| {\Gamma}|}},
\end{equation}
therefore $d\log \hz=0$ by (\ref{tors}). 
$\Box$

At this point, it is perhaps appropriate to mention the
behaviour of the invariant under the two operations of 
reversing the orientation and taking the mirror image of the knot. 
If $\phi \in \cal K$ is a representative of a knot, denote by 
$\phi^* =-\phi$ the mirror image of this knot 
and $\bar\phi(s) =\phi (1-s)$ the knot with the opposite 
orientation.
Let $\bar\Gamma$ be the graph obtained from $\Gamma$ by reversing the 
orientation of the Wilson line.
We have the following properties:
\bea
I(\Gamma) (\phi^*) & = & (-1)^{{\mathrm deg} \Gamma}I(\Gamma) (\phi),\\
I(\Gamma) (\bar\phi) & = & (-1)^{n}I(\bar\Gamma)(\phi).
\eea
The first property tells us that an invariant of even order 
appearing in the Chern-Simons expansion cannot distinguish a knot from 
its mirror image. This property can also be recasted in the usual form:
\begin{equation}
Z(\phi^*) (\hbar) = Z(\phi) (-\hbar).
\end{equation}
From the second property we can deduce that, 
if the space ${\cal A}$ of BN diagrams is such that 
\begin{equation}
D(\Gamma) = (-1)^n D(\bar\Gamma),
\label{even}
\end{equation}
then the invariant $\hz$ does not depend on the orientation of knots.
It is well-known \cite{bn2} that this holds in ${\cal A}^*$ for all weight 
systems constructed from simple Lie algebras, but, as 
far as the authors are aware, whether (\ref{even}) is true or not
is still an open question.
All these remarks also apply to the universal invariant 
constructed by Kontsevich.

\begin{figure}[htbp]
$$\epsf{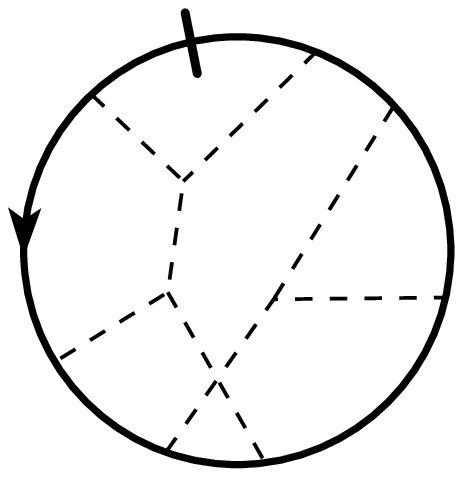}[xscale=2/3,yscale=2/3] $$ 
\caption{\label{gmark}}
\end{figure}

\begin{figure}[htbp]
$$\epsf{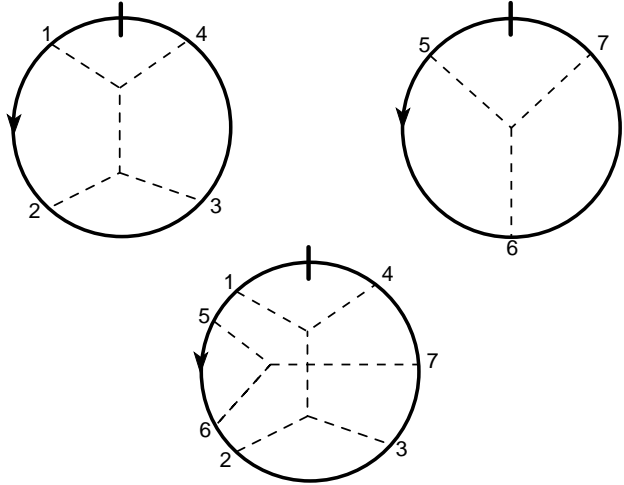}[xscale=2/3,yscale=2/3]$$ 
\caption{$\sigma(1,2,3,4,5,6,7)=(1,4,5,7,2,3,6)$\label{gshuf}}
\end{figure}

\begin{figure}[htbp]
$$\epsf{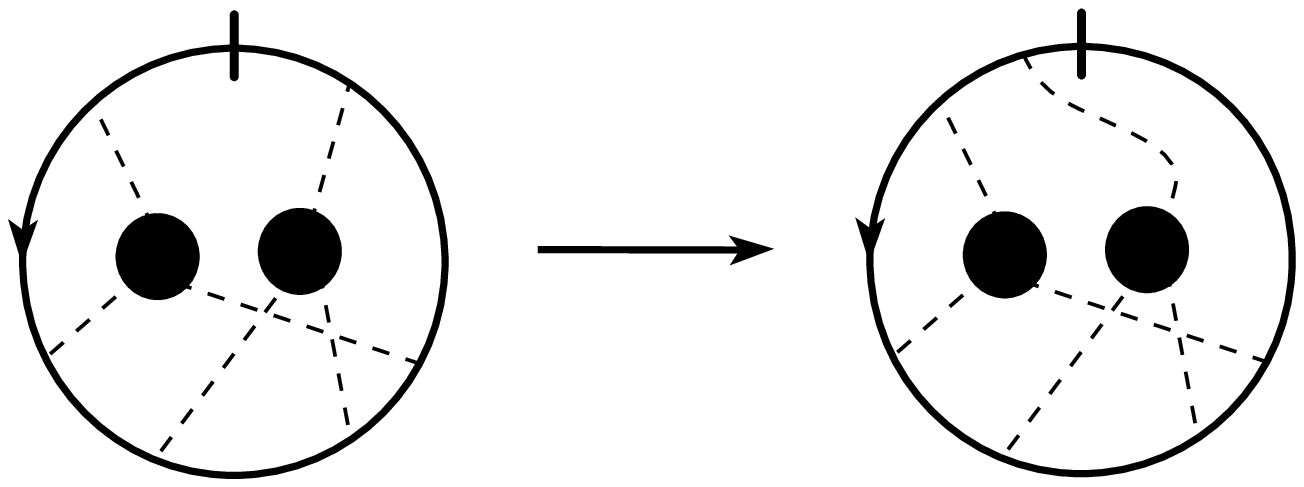}[xscale=2/3,yscale=2/3]$$ 
\caption{\label{figtau}}
\end{figure}

The remainder of this section is devoted to the
proof of lemma \ref{grlike}, which requires some preparation.
The main ideas go back to section 9.6 of Bar-Natan's thesis \cite{bn1}.
First, we need to define {\em marked graphs}. They are pairs
$(\Gamma, e)$, where $\Gamma$ is a trivalent graph and $e\in E^W_\Gamma$
is an edge of $W_\Gamma$. The edge $e$ is also called the marking
of $(\Gamma, e)$. Marked graphs are depicted as in fig. \ref{gmark}.
Recall that for trivalent graphs $\Gamma$, $I(\Gamma)$ is an integral
over the configuration space $C_{n,t}$, where points on the circle 
are cyclically ordered. For marked graphs $(\Gamma, e)$, we define 
$I(\Gamma,e)$ to be the same integral, but taken over a configuration
space $L_{n,t}$, where we consider the points on the circle 
linearly ordered ($0<s_1<\cdots<s_n<1$) and the marked Wilson edge 
corresponds to the interval $[s_n,s_1]$. Thus, by definition we have:
\begin{equation}
I(\Gamma) = \sum_{e\in E^W_\Gamma} I(\Gamma,e).
\end{equation}
The group ${\mathrm Aut}(\Gamma)$ acts on $E^W_\Gamma$, and
if $g\in {\mathrm Aut}(\Gamma)$, $I(\Gamma, g\cdot e)=I(\Gamma,e)$.
Let ${\mathrm Aut}(\Gamma)_e$ be the stabilizer of $e\in E^W_\Gamma$,
and $\Gamma' = \Gamma - W_\Gamma$ be the interior of $\Gamma$.
It is easy to check that
\begin{equation}
{\mathrm Aut}(\Gamma)_e = {\mathrm Aut}(\Gamma)_{e'} = {\mathrm Aut}(\Gamma'),
\label{aute}
\end{equation}
for all $e, e'\in E^W_\Gamma$. Let
$M(\Gamma)=E^W_\Gamma/{\mathrm Aut}(\Gamma)$ be the set of 
equivalence classes of markings of $\Gamma$. The number of markings
in the class of any $e\in E^W_\Gamma$ is 
$|{\mathrm Aut}(\Gamma)/{\mathrm Aut}(\Gamma')|= |\Gamma|/|\Gamma'|$.
Therefore,
\begin{equation}
\frac{I(\Gamma)}{|\Gamma|} = \sum_{e\in M(\Gamma)} 
\frac{I(\Gamma,e)}{|\Gamma'|}.
\label{sumint}
\end{equation}
Given two marked graphs $(\Gamma_1,e_1)$ and $(\Gamma_2,e_2)$,
with $n_1$ and $n_2$ external vertices, and a $(n_1,n_2)$\/-shuffle $\sigma$,
we can construct a new marked graph 
$(\Gamma,e) = (\Gamma_1,e_1) \shuf (\Gamma_2,e_2)$, which we call
a shuffle product. A $(n_1,n_2)$\/-shuffle is a permutation $\sigma$
of $\{1,\ldots,n_1+n_2\}$ such that  $\sigma(1)<\cdots <\sigma(n_1)$
and $\sigma(n_1+1)<\cdots <\sigma(n_1+n_2)$. The set of 
$(n_1,n_2)$\/-shuffles will be denoted by $\Sigma(n_1,n_2)$.
The shuffle product is defined as follows: label the external
vertices of $(\Gamma_1,e_1)$ by $1,2,\ldots,n_1$, going around
the Wilson line in the sense given by its orientation, so that
for the oriented edge $e_1$, $\partial e_1 = (n_1,1)$, for the next
edge $e'_1$, $\partial e'_1 = (1,2)$, etc. Similarly, label
the external vertices of $(\Gamma_2,e_2)$ by $n_1+1,\ldots,n_1+n_2$,
with $\partial e_2 = (n_1+n_2, n_1+1)$. The interior of the
shuffle product is $\Gamma' = \Gamma'_1 \cup \Gamma'_2$.
Going around the Wilson line, its external vertices are 
$\sigma^{-1}(1),\sigma^{-1}(2),\ldots,\sigma^{-1}(n_1+n_2)$,
and the marked edge is $\partial e = (\sigma^{-1}(n_1+n_2),\sigma^{-1}(1))$.
This is illustrated in fig. \ref{gshuf} .
Now it should be fairly obvious that
\begin{equation}
I(\Gamma_1,e_1) I(\Gamma_2,e_2) = \sum_{\sigma\in\Sigma(n_1,n_2)}
I((\Gamma_1,e_1) \shuf (\Gamma_2,e_2)).
\label{prodint}
\end{equation}
Let 
$S(\Gamma_1,\Gamma_2)=M(\Gamma_1)\times M(\Gamma_2)\times \Sigma(n_1,n_2)$.
We have a map 
$s : S(\Gamma_1,\Gamma_2)\ni (e_1,e_2,\sigma) \mapsto 
(\Gamma_1,e_1) \shuf (\Gamma_2,e_2)$. Let
\begin{equation}
P(\Gamma|\Gamma_1,\Gamma_2)=\{t\in S(\Gamma_1,\Gamma_2) \, | \, 
s(t)=(\Gamma,e), \mbox{ for some marking }e\},
\end{equation}
and $n(\Gamma|\Gamma_1,\Gamma_2)=|P(\Gamma|\Gamma_1,\Gamma_2)|$.
Elements of $P(\Gamma|\Gamma_1,\Gamma_2)$ will be called
partitions of $\Gamma$ in two parts $\Gamma_1,\Gamma_2$.
Clearly, partitions can be defined for an arbitrary number
of parts.
\begin{lemma}
Let $\Gamma_1, \ldots , \Gamma_m \in G^3$, then 
\begin{equation}
{I(\Gamma_1)\over|\Gamma_1|} \cdots {I(\Gamma_m)\over|\Gamma_m|} =
\sum_{\Gamma \in G^3}
 {I(\Gamma)\over |\Gamma|} \,
n(\Gamma|\Gamma_1, \ldots , \Gamma_m).
\end{equation}
\label{prodlem}
\end{lemma}
\noindent
{\em Proof.} We can restrict ourselves to the case $m=2$.
Using (\ref{aute}), (\ref{sumint}) and (\ref{prodint}) we get
\begin{equation}
{I(\Gamma_1)\over|\Gamma_1|} {I(\Gamma_2)\over|\Gamma_2|} =
\frac{1}{|\Gamma'_1| |\Gamma'_2|} 
\sum_{e_1\in M(\Gamma_1)} \, \sum_{e_2\in M(\Gamma_2)} \,
\sum_{\sigma\in\Sigma(n_1,n_2)} I((\Gamma_1,e_1) \shuf (\Gamma_2,e_2)).
\label{prod2int}
\end{equation}
Observe that if $(\Gamma,e) = (\Gamma_1,e_1) \shuf (\Gamma_2,e_2)$,
${\mathrm Aut}(\Gamma')={\mathrm Aut}(\Gamma'_1)\times
{\mathrm Aut}(\Gamma'_2)$, so that 
\begin{equation}
|\Gamma'_1| |\Gamma'_2| = |\Gamma'|.
\label{remord}
\end{equation}
Consider the map
$\tau : S(\Gamma_1,\Gamma_2) \rightarrow S(\Gamma_1,\Gamma_2)$
defined by fig. \ref{figtau}. 
It is invertible and generates
a finite cyclic group $G$ of order at most $n_1+n_2$. There is
also a transformation, which we denote by the same letter $\tau$,
acting directly on marked graphs: $\tau\cdot (\Gamma,e)=(\Gamma,e')$,
where $e'$ is the marking directly adjacent to $e$.
Moreover, $\tau\circ s = s\circ\tau$.

The sum on the r.h.s. of (\ref{prod2int}) can be rearranged into
a sum over the orbits of $G$ in $S(\Gamma_1,\Gamma_2)$.
If $(e_1,e_2,\sigma)\in S(\Gamma_1,\Gamma_2)$, 
$s((e_1,e_2,\sigma))=(\Gamma,e)$, then the orbit $G\cdot(e_1,e_2,\sigma)$
covers a certain number of times, say
$n(e_1,e_2,\sigma)$, the orbit $G\cdot(\Gamma,e)$.
It follows that the contribution of one orbit is:
\begin{equation}
\sum_{(f_1,f_2,\mu)\in G\cdot(e_1,e_2,\sigma)} 
I((\Gamma_1,f_1) \times_\mu (\Gamma_2,f_2)) =
n(e_1,e_2,\sigma) \sum_{e\in M(\Gamma)} I(\Gamma,e).
\end{equation}
Finally, if we add
the contributions of all the orbits in $S(\Gamma_1,\Gamma_2)$
which cover the same $G\cdot(\Gamma,e)$, we get
\begin{equation}
n(\Gamma|\Gamma_1,\Gamma_2) \sum_{e\in M(\Gamma)} I(\Gamma,e),
\end{equation}
where $n(\Gamma|\Gamma_1,\Gamma_2)$ is the sum of all the
factors $n(e_1,e_2,\sigma)$. Taking (\ref{remord}) into
account and using (\ref{sumint}) once again, the proof is completed.
$\Box$

\noindent
{\em Proof of lemma \ref{grlike}.}
Writing
\begin{equation}
Z = 1 + \sum_{\Gamma \in G^3_*} \hbar^{{\mathrm deg}\Gamma} 
{I(\Gamma)\over |\Gamma|} D(\Gamma),
\end{equation}
and using the product in $\cal A$, we define $\log Z$ by the
formal power series expansion
\begin{equation}
\log Z = \sum_{m=1}^{\infty} \frac{(-1)^{m+1}}{m}
\left( 
\sum_{\Gamma \in G^3_*} \hbar^{{\mathrm deg}\Gamma} 
{I(\Gamma)\over |\Gamma|} D(\Gamma)
\right)^m .
\end{equation}
Now from lemma \ref{prodlem} we get
\begin{eqnarray}
\left( 
\sum_{\Gamma \in G^3_*} \hbar^{{\mathrm deg}\Gamma} 
{I(\Gamma)\over |\Gamma|} D(\Gamma)
\right)^m & = & \sum_{\Gamma_1,\ldots,\Gamma_m}
{I(\Gamma_1)\over|\Gamma_1|} \cdots {I(\Gamma_n)\over|\Gamma_m|} \,
\prod_{j=1}^m D(\Gamma_j) \hbar^{{\mathrm deg}\Gamma_j}
\nonumber \\
 & = & \sum_{\Gamma} \hbar^{{\mathrm deg}\Gamma}\frac{I(\Gamma)}{|\Gamma|}
       \sum_{\Gamma_1,\ldots,\Gamma_m} n(\Gamma| \Gamma_1,\ldots,\Gamma_m)
       \prod_{j=1}^m D(\Gamma_j).
\end{eqnarray}
Thus, $\log Z$ becomes
\begin{equation}
\log Z = \sum_{\Gamma \in G^3_*} \hbar^{{\mathrm deg}\Gamma} 
\frac{I(\Gamma)}{|\Gamma|} \, c(\Gamma),
\end{equation}
where
\begin{equation} 
c(\Gamma) = \sum_{\mbox{\scriptsize partitions } P} 
\frac{(-1)^{|P|+1}}{|P|} \, D(P),
\label{spart}
\end{equation}
the sum is over all partitions $P$ of $\Gamma$, $|P|$ is the number
of parts, and $D(P)=\prod_{j=1}^m D(\Gamma_j)$ if $P$ is a partition
with $|P|=m$ parts $\Gamma_1,\ldots,\Gamma_m$. It remains to 
show that $c(\Gamma)=C(D(\Gamma))$. It is clear that 
if $\Gamma$ is primitive, $c(\Gamma)=D(\Gamma)$.
If we can prove that $c(\Gamma)=0$ when $\Gamma$ is not prime,
then by lemma \ref{lemC} we are done. Bar-Natan has 
shown in \cite{bn1} how to achieve this in three steps, which
we reproduce here since his thesis has not been
published. 

\noindent 
{\em Step 1.} We have 
$n(\Gamma|\Gamma_1,\Gamma_2)=n(\Gamma|\Gamma_2,\Gamma_1)$, and
more generally
$n(\Gamma|\Gamma_1,\Gamma_2,\ldots,\Gamma_m)=
n(\Gamma|\Gamma_m,\Gamma_1,\ldots,\Gamma_{m-1})$. Define a cyclic
partition of $\Gamma$ to be a partition modulo a cylic permutation
of the parts. Then the total number of partitions into $m$ parts is
divisible by $m$ and is the number of cyclic partitions into $m$ parts.
Thus we can rewrite (\ref{spart}) as
\begin{equation} 
c(\Gamma) = \sum_{\mbox{\scriptsize cyclic partitions } P} 
(-1)^{|P|+1} \, D(P).
\end{equation} 

\noindent
{\em Step 2.} We show that $c(\Gamma)=0$ if
$\Gamma$ is the connected sum along the Wilson lines of two graphs
$\Lambda$ and $M$, such that $\Lambda$ is primitive. 
To do this, we construct
an involution $\rho$ on the set of cyclic partitions
of $\Gamma$ such that $|\rho(P)|=|P|\pm 1$ and $D(\rho(P))=D(P)$.
Take a cyclic partition $P=\{\Gamma_1,\ldots,\Gamma_m\}$.
We can always assume that $\Lambda\subset\Gamma_1$. 
If $\Lambda=\Gamma_1$, we put 
$\rho(P)=\{\Lambda\cup\Gamma_2,\ldots,\Gamma_m\}$,
otherwise $\Lambda$ is properly contained in $\Gamma_1$ and we put
$\rho(P)=\{\Lambda,\Gamma_1-\Lambda,\Gamma_2,\ldots,\Gamma_m\}$.

\noindent
{\em Step 3.} We show that $c(S)=c(T)-c(U)$, where $S,T,U$ are the
three graphs appearing in the STU relation of \cite{bn2}.
The argument is the following: 
to every partition of $S$ there corresponds a partition of $T$ and $U$,
but the converse is not true. However, the only partitions of $T$ and $U$
to which one cannot associate a partition of $S$ are those where
the two edges not contained in the Wilson line belong to two
distinct parts. Call these the exceptional partitions. There is
an obvious 1-1 correspondence between the exceptional partitions
of $T$ and $U$, which preserves the number of parts. Hence the
contributions of the exceptional partitions cancel in $c(T)-c(U)$.

Now by step 3, there is a map $\tilde{c} : {\cal A} \rightarrow {\cal A}$
such that $c(\Gamma)=\tilde{c}(D(\Gamma))$, for all $\Gamma\in G^3$.
If $\Gamma$ is not prime, $D(\Gamma)=D(\Lambda)D(M)$, with 
both factors non-trivial. But in ${\cal A}$, $D(\Gamma)$ is
a linear combination of products of primitive elements, therefore
$c(\Gamma)=0$ by step 2.
$\Box$.

\section{$\hz$ is a universal Vassiliev invariant}
In this section we show that $\hz$ is a universal Vassiliev Invariant.
First of all, let us briefly recall what this means.
Put $\hz(K)=\hz(\phi,\nu)$, where $K$ is a framed knot, and
\begin{equation}
\hz(K) = \sum_{N\geq 0} \hz_N(K) \hbar^N.
\end{equation}
We extend $\hz(K)$ to singular knots in the usual
way: if $K^j$ is a singular knot with $j$ double points,
we define $\hz(K^j)$ by
\begin{equation}
\hz_N(K^j)=\hz_N(K^{j-1}_+)-\hz_N(K^{j-1}_-),
\label{dsing}
\end{equation}
\begin{figure}[htbp]
$$\epsf{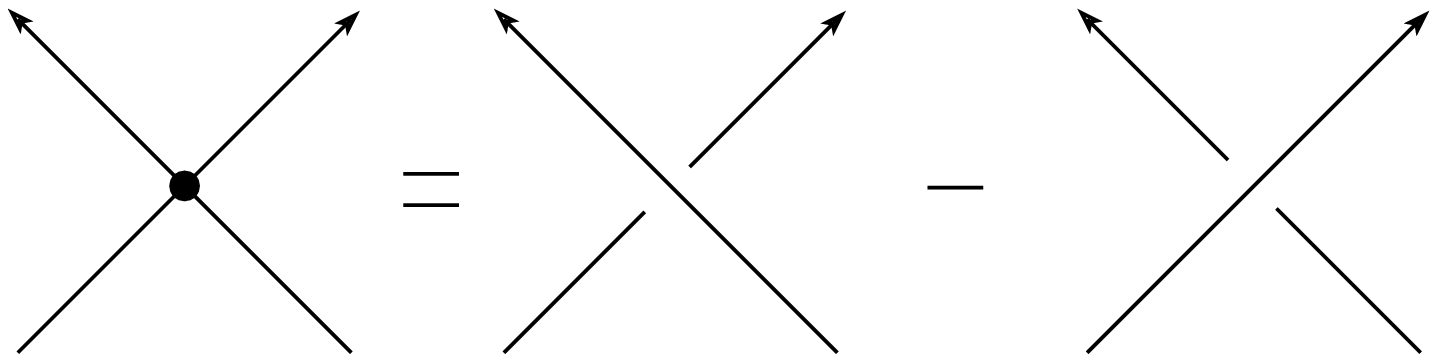}[xscale=1/2,yscale=1/2] 
$$ \caption{\label{singul}}
\end{figure}
where $K^{j-1}_\pm$ are two knots with $j-1$ double points
obtained by desingularizing one of the double points of $K^j$
as shown in fig. \ref{singul}. We also need the fact that
each singular knot $K^j$ defines a unique chord diagram $\Gamma(K^j)$
of degree $j$ \cite{bn2}. A formal power series 
$\Phi(K)=\sum_{N\geq 0} \Phi_N(K) \hbar^N \, \in {\cal A}[[\hbar]]$ 
is a universal Vassiliev invariant if it satisfies the
following two properties:
\begin{description}
\item[U1] $\Phi_N(K^j)=0$ if $j>N$, i.e. $\Phi_N(K)$ is a Vassiliev
invariant of degree $\leq N$,
\item[U2] $\Phi_N(K^N)=D(\Gamma(K^N))$, i.e. 
$\Phi(K^N)-\hbar^N D(\Gamma(K^N))$ is divisible by $\hbar^{N+1}$.
\end{description}
\begin{theo}
\label{universal}
$\hz(K)$ is a universal Vassiliev invariant.
\end{theo}
{\em Proof.}
If we desingularize the double points of $K^j$ in all possible ways,
we obtain $2^j$ knots $K_\choice$,  where
$\choice=(\varepsilon_1,\ldots,\varepsilon_j)$ and
$\varepsilon_i=\pm 1$ for all  $1\leq i \leq j$. It is possible
to choose the corresponding embeddings and framings $\phi_\choice$
and $\nu_\choice$ such that the total torsion 
$\tau(\phi_\choice,\nu_\choice)$ is independent of $\choice$.
Thus in the forthcoming arguments, we can ignore the factor
containing the torsion in $\hz(K)$.

For non-singular $K$, $\hz_N(K)$  is a sum of
contributions $I(K,\Gamma)$ indexed by graphs $\Gamma$
of degree $N=(n+t)/2$, where $n$ is the number of (external)
vertices on the Wilson line $W_\Gamma$ and $t$ is the number of
(internal) trivalent vertices. 

By (\ref{dsing}), we have for each $\Gamma$ a contribution 
$I(K^j,\Gamma)$ to $\hz_N(K^j)$,
which is an alternating sum of
$2^j$ terms $I(K_\choice,\Gamma)$, corresponding to all possible 
desingularizations of $K^j$.
Let $\choice(i)$ denote the vector obtained from $\choice$
by changing the sign of the $i$\/-th coordinate. Then for $1\leq k\leq j$,
\begin{eqnarray}
I(K^j,\Gamma) & = & \sum_\choice (\prod_{i=1}^j \varepsilon_i)
 I(K_\choice,\Gamma) \nonumber\\
 & = & \sum_{\choice, \:\varepsilon_k=1} (\prod_{i\neq k} \varepsilon_i)
 (I(K_\choice,\Gamma)-I(K_{\choice(k)},\Gamma)) \label{expand}.
\end{eqnarray}
Now $I(K_\choice,\Gamma)$ is the integral of a function
of $n$ real variables
$s_\alpha$, $1\leq\alpha\leq n$, which are parameters
along the knot $K_\choice$, and $t$ points of $\R^3$. 
Let $R_i$, $1\leq i \leq j$ be a
small ball containing the $i$\/-th singular point of $K^j$.
Let $C_{\choice,i}=R_i\cap K_\choice$. 
By definition,
\begin{equation} 
K_\choice-C_{\choice,i}=K_{\choice(i)}-C_{\choice(i),i}, 
\label{toto}
\end{equation}
for each $i$.
Let
\begin{equation}
D_i=\{ s\in[0,1] \, | \, K^j(s)\in R_i\}.
\end{equation}
Remembering the discussion surrounding (\ref{plonge}), we can identify 
$D_i$ with a subset of $W_\Gamma$.
Let us say that $D_i$ is {\em occupied} if it contains at least one
of the $n$ external vertices. 
Let $U_\choice$ be the integration domain of $I(K_\choice,\Gamma)$,
$E_k=\{ s_\alpha\notin D_k, \;\forall\alpha\}$. If $F_k\subset E_k$,
then $U_\choice \cap F_k$ is a sub-domain with $D_k$ not occupied.
Denote by $I(K_\choice,\Gamma,F_k)$
the contribution of $U_\choice \cap F_k$ to $I(K_\choice,\Gamma)$.
We can also define $I_k(K^j,\Gamma,F_k)$, replacing $I(\cdot,\Gamma)$ 
by $I(\cdot,\Gamma,F_k)$
everywhere in the first equality in (\ref{expand}).
\begin{lemma}
If $j\geq k\geq 1$ and $F_k\subset E_k$, $I(K^j,\Gamma,F_k)=0$.
\end{lemma}
{\em Proof.} By (\ref{toto})
and the second equality of (\ref{expand}), the lemma
becomes obvious. $\Box$

This lemma already implies that $\hz_N(K)$ is a finite type
invariant. Indeed, it says that all non-zero contributions to
$I(K^j,\Gamma)$ come from the domain ${\cal D}_\Gamma$
where all the $j$ regions $D_i$ are occupied. 
But if $\Gamma$ is such that $j>n$, then ${\cal D}_\Gamma=\emptyset$
and $I(K^j,\Gamma)=0$. Since for all $\Gamma$ of degree $N$, $2N\geq n$,
we get $\hz_N(K^j)=0$ if $j>2N$. We need to work a little bit more
to prove the theorem. 

Put $\bar{D}=W_\Gamma-\bigcup_i D_i$. It corresponds
to the common part of all the $K_\choice$.
Let $\Gamma'=\Gamma-\{ {\mathrm edges\;\;of\;\;} W_\Gamma\}$ 
be the interior of $\Gamma$. 
\begin{lemma}
If $\Gamma$ is a graph of degree $N$ with $j>N$, 
and every $D_i$ is occupied, then
there exist $D_i$ and $D_m$, $i\neq m$, such that there is
a path in $\Gamma'$ going from $D_i$ to $D_m$.
\label{lem1}
\end{lemma}
{\em Proof.} 
Let $\Gamma_i$ be the 
connected component of $\Gamma'$ containing all vertices from $D_i$.
If there is no pair $(i,m)$, $i\neq m$ with a path from $D_i$ to $D_m$,
then for all $i$, the external vertices of $\Gamma_i$ all
belong to $D_i$ or $\bar{D}$, but never to $D_m$ with $m\neq i$.
Thus the graphs $\Gamma_i$ are all disjoint. Since every $\Gamma_i$ has
at least two vertices, and there are $j$ such graphs,
$N\geq (2j)/2=j$. $\Box$

For a while, let us forget knots and consider embeddings in $\R^3$
of connected graphs $\Gamma$, whose vertices are either trivalent (internal)
or univalent (external). Assign fixed locations $x_\alpha\in \R^3$,
$\alpha=1,\ldots,n$ to the external
vertices of $\Gamma$. Integrating over all the internal vertices
the form $\omega(\Gamma,\orga)$, we obtain a real-valued
$n$\/-form $g_{\,\Gamma}(x)$, where $x=(x_1,\ldots,x_n)$.
\begin{lemma}
Let $x_{\alpha\beta}=x_\alpha-x_\beta$. If $\alpha\neq\beta$, 
\[
\lim_{{|x_{\alpha\beta}|}\to\infty} g_{\,\Gamma}(x)=0.
\]
\label{lem2}
\end{lemma}
{\em Proof.}
Consider a path in $\Gamma$ connecting the two vertices $x_\alpha$
and $x_\beta$. 
If this path consists of a single edge, the lemma
is trivial, so we assume that it passes through internal vertices
$z_1,z_2,\ldots,z_k$, with $k\geq 1$, which are numbered in such 
a way that $z_j$ is the vertex reached after traveling through
$j$ consecutive edges, starting from $x_\alpha$. 

Now make the
substitution $z_j \rightarrow z_j+x_\beta$ of the integration
variables. Then in $g_{\,\Gamma}$ the $z_1$ integral is
$g_Y(x_{\alpha\beta},z_2,w)$, where $w$ is the end of the third
edge connected to $z_1$, and $g_Y$ is the form corresponding
to the graph $Y$. Fortunately, $g_Y$ was computed explicitly
in \cite{guad}, and looking at the expression one can
see easily that $g_Y(x_{\alpha\beta},z_2,w)\rightarrow 0$
as ${|x_{\alpha\beta}|}\to\infty$. $\Box$

U1 now follows from the last two lemmas: 
by the invariance
of $\hz_N(K)$, we can assume that the balls $R_i$ are all very far
from each other. Then by lemma \ref{lem1},
each term $I(K_\choice,\Gamma)$ contributing
to $\hz_N(K^j)$ will contain a factor $g_\Delta(x)$, where $\Delta$
is a connected subgraph of $\Gamma'$ having two external vertices
which are very far apart, and by lemma \ref{lem2} it vanishes.

Let us now prove U2. We use the
same notations as in the proof of lemma \ref{lem1}. 
If there is a path in $\Gamma'$ connecting disjoint regions
$D_i$, $D_m$, $m\neq i$, then by lemma \ref{lem2} $I(K^j,\Gamma)=0$.
Thus the only graphs $\Gamma$ which contribute to $\hz_N(K^j)$ are
those such that all the subgraphs $\Gamma_i$ are disjoint. 
Since $N=j$, every $\Gamma_i$ has exactly two vertices.

Using the invariance of $\hz_N(K)$, we can assume that
for each $i$, there is a very large neighborhood $V_i$ of $R_i$ in
which the two connected components of 
$L_{\choice,i}=K_\choice\cap V_i$ are both
planar curves, for all $\choice$. 
Then if for some value of $i$,
the unique edge of $\Gamma_i$ joins two points of the same component 
of $L_{\choice,i}$, $I(K^j,\Gamma)=0$, by planarity.
Therefore all $\Gamma_i$ join two distinct connected components
of $L_{\choice,i}$. Thus we see that the only non-vanishing
contribution to $\hz_N(K^j)$ for $N=j$ comes from the graph
$\Gamma(K^j)$. 
It is easy to see that each $V_i$ will contribute a factor 
$1$ to a product decomposition of $I(K^j,\Gamma(K^j))$
corresponding to the $j$ pairs of integration variables $s_\alpha$.
This concludes the proof of the theorem. $\Box$

{\bf Acknowledgements.} We would like to thank Fran\c{c}ois Delduc and
Raymond Stora for many useful discussions, and the referee
for his helpful suggestions and comments.
D.A. thanks ENSLAPP  and L.F. thanks ETH Zurich for their hospitality. 

\setcounter{section}{0}
\renewcommand{\thesection}{Appendix \Alph{section}}
\renewcommand{\theequation}{\Alph{section}.\arabic{equation}}

\section{Integration along the fiber}
\label{appa}
Let $(M,B,p)$ be a bundle ($p : M \rightarrow B$) such that  the fiber 
$p^{-1}(x)=F_x,\, x\in B, $ is a compact finite-dimensional 
oriented (possibly stratified) manifold.
The push-forward, or integration along the fiber, is a linear morphism
$p_*^M$ from the space of $(n+p)$\/-differential forms on $M$, 
to the space of $p$\/-differential forms on $B$, where $n$ is the 
dimension of the fiber. 
It is defined as follows:
let $x\in B$ and  $X_1,\cdots,X_p$  be $p$ vectors in $ T_xB$, 
let $\omega $ be a $(n+p)$\/-form on $M$, then the push-forward is 
\begin{equation}
(p_*^M \omega)_x (X_1,\cdots,X_p) =
\int_{F_x} 
i_{\tilde{X}_p}\cdots i_{\tilde{X}_1}
\omega,
\end{equation}
where $i$ denotes the interior product, and $\tilde{X}_i$ 
 is any lift of $X_i$ as a section of $TM$ over $F_x$.
This definition is independent of the lift's choice.
This integration along the fiber is a direct generalization of 
usual integration of forms on compact manifolds:
if $\omega$ is a $n$\/-form then the push-forward is
\begin{equation}
(p_*^M \omega) =\int_{F_x} 
\omega,
\end{equation}
and $B$ appears in this case as a parameter space.
An important property of the push-forward is its commutation with the 
differential operator, leading to the generalization of Stokes theorem:
\begin{equation}
d_B p_*^M \omega= p_*^M d_M \omega+ 
(-1)^{{\mathrm deg}( p_*^{\partial M}\omega)} p_*^{\partial M}\omega.
\label{diff}
\end{equation}
Here $d_B$ (resp. $d_M$) denotes the differential operator on $B$ (resp. $M$).
In order to define the integration $p_*^{\partial M}$, we need to 
provide the boundary of the fiber $\partial {F_x}$ with an orientation
(in the case of a stratified space the boundary is the codimension 
$1$ strata). This orientation is induced by the orientation on the whole fiber 
as follows:
let $\Omega $ be an orientation (i.e. an $n$\/-form) on the fiber 
${F_x}$, and $n_e$ the in-going normal vector field of the 
boundary. The boundary orientation is defined by $i_{n_e} \Omega$.
The formula (\ref{diff}) holds with this convention.
Let us remark that it is not the usual orientation convention of the
boundary, for if we apply (\ref{diff}) to an $n$\/-form $\omega$ on a trivial
bundle $M=B\times F$ 
which does not depend on the parameter space $B$, then the l.h.s. of 
(\ref{diff}) is zero, and this equation reduces to Stokes' formula 
\begin{equation}
\int_F d\omega =-\int_{\partial F} \omega.
\end{equation}
Another important property of the push-forward is its behaviour with 
respect to the pullback. If $\phi$ is a form on $B$, $\omega$ a form on $M$, 
and $p^*(\phi)$ denotes the pullback of $\phi$ on $M$, then 
\begin{equation}
p_*^M(p^*(\phi)\wedge \omega)=\phi \wedge p_*^M(\omega).
\end{equation}

\end{document}